\documentclass[11pt,twoside]{article}
\usepackage{latexsym}
\usepackage{geometry}                
\usepackage{color}
\usepackage[english]{babel}

\usepackage{amsmath}
\usepackage{subeqnarray}
\usepackage{graphicx}
\usepackage{amssymb}
\usepackage{epstopdf}

\usepackage{amsfonts}
\usepackage{cases}
\newcommand{\commentout}[1]{}
\newcommand{\R}{\mathbb{R}}

\newcommand {\eps}  {\varepsilon}

\newcommand {\Chi} {{\bf \raise 2pt \hbox{$\chi$}} }
\newcommand {\f}   {\frac}
\newcommand {\p}   {\partial}
\newcommand{\ds}{\displaystyle}
\newcommand{\ud}{\, \mathrm{d}}

\newcommand{\qed}{{ \hfill 
{\unskip\kern 6pt\penalty 500 
\raise -2pt\hbox{\vrule\vbox to 6pt{\hrule width 6pt
\vfill\hrule}\vrule} \par} }}
\newtheorem{Def}{\indent{Definition}}[section]
\newtheorem{prop}{\indent{Proposition}}[section]
\newtheorem{theor}{\indent{Theorem}}[section]

\newtheorem{lem}{\indent{Lemma}}[section]
\newtheorem{Rmq}{\indent{Remark}} 
\title{ A General Inverse Problem for the Growth-Fragmentation Equation}
 \author{M. Doumic \thanks{INRIA Paris-Rocquencourt, EPI BANG, 
           Domaine de Voluceau, F~78153 Le Chesnay cedex;
email: marie.doumic@inria.fr} \and L\'{e}on M. Tine\footnotemark[2] \footnotemark[3] \thanks{Laboratoire d'Analyse Num\'erique et d'Informatique (LANI) Universit\'e Gaston Berger, B.P. 234 Saint-Louis, S\'en\'egal}\thanks{Labo P. Painlev\'e UMR 8524 CNRS \& Universit\'e des Sciences et Technologies Lille1}\thanks{Project-Team SIMPAF, INRIA Lille Nord Europe Research Center, Park Plazza, 40 avenue Halley, 59650 Villeneuve d'Ascq cedex}}
\begin{document}
\maketitle
\begin{abstract}
Aggregation-fragmentation equations arise in many different contexts, ranging from cell division, protein polymerization, biopolymers, neurosciences etc. Direct observation of temporal dynamics being often difficult, it is of main interest to develop theoretical and numerical methods to recover reaction rates and parameters of the equation from indirect observation of the solution. Following the work done in \cite{P.Z} and \cite{D.P.Z} for the specific case of the cell division equation, we address here the general question of recovering the fragmentation rate of the equation from the observation of the time-asymptotic solution, when the fragmentation kernel and the growth rates are fully general. We give both theoretical results and numerical methods, and discuss the remaining issues.
\end{abstract}
%
  
{\bf keywords.} {Growth-Fragmentation equation ; general fragmentation kernels ; inverse problem ; eigenvalue problem.}  

{\bf MSC Classification.} 35Q92, 35R30, 45Q05.

\section*{Introduction}
To model the behavior of a population where growth and division depend on a  structuring quantity of the individuals such as size, the following mass-balance equation is currently used:
\begin{equation} \label{P_t}
\left\{\begin{array}{ll}
 \ds \frac{\partial}{\partial_{t}}n(t,x)+c\frac{\partial}{\partial_{x}}(g(x)n(t,x))+B(x)n(t,x)=2\int_{0}^{+\infty}B(y)\kappa(x,y)n(t,y)\ud y, \; \;  t >0, \;  x>0, \\
 \\
  n(t=0,x)=n^{0}(x), \quad  x\geq0, \\
\\
g(x=0)  n(t,x=0)=0, \, t\geq0, \;  \; c>0. \\
   \end{array} \right.\vspace{0.1cm}
\end{equation}
Here, $n$ denotes the density of the individuals structured by the size variable $x$ at time $t;$ the growth rate is given by $g(x)$; the division rate $B(y) \kappa(x,y)$ represents the rate at which a given individual of size $y$ gives birth to two individuals of size respectively $x$ and $y-x,$ whereas $B(y)$ is the total rate of division for individuals of size $y.$ This physical interpretation of $\kappa(x,y)$ leads to the following assumptions:
\begin{equation}\label{as:kappa}
\kappa(x,y)=0 \;\;  \forall x>y, \qquad \ds\int_{0}^{\infty} \kappa(x,y)\ud x =1,\qquad \kappa(x,y) = \kappa(y-x,y).
\end{equation}
By simple integration and symmetry, it leads to the following well-known relation, expressing the conservation of mass by the division process:
$$ 
\int_{0}^{\infty} x\kappa(x,y)\ud x =\frac{y}{2}.$$
Problem \eqref{P_t} or its variants arises in many different contexts, ranging from cell division, protein polymerization, telecommunication, neurosciences, and its mathematical study can provide useful information on the qualitative behavior of the phenomenon under consideration (see, among many others, \cite{MetzDiekmann,P}). To be able to use it as a predictive model however, it is crucial to be able to estimate quantitatively its parameters $g,$ $B$ and $\kappa.$  

A first step consists in the use of the asymptotic behavior of this equation, as first proposed in \cite{P.Z}. Indeed, by general relative entropy principle it is proven (see e.g \cite{P_ryz, M.S.P, D.G}) that under suitable assumptions on $\kappa$, $g$ and $B$ one has
$$\int\limits_0^\infty \vert n(t,x) e^{-\lambda_0 t} - \rho_0 N(x) \vert \phi(x) dx \underset{t\to \infty}{\longrightarrow} 0,$$
with $\rho_0=\int n^0(x) \phi(x) dx$ and $(\lambda_{0}, N, \phi)$ is the unique eigenpair solution of the following problem:
 \begin{equation}\label{(P)}
\left\{\begin{array}{ll}
  c\ds\frac{\partial}{\partial_{x}}(g(x)N(x))+(B(x)+\lambda_{0})N(x)=2\ds\int^{+\infty}_{0}B(y)\kappa(x,y)N(y)\ud y, \\
  \vspace{0.01cm}
  gN(x=0)=0; \;  N(x)\geq0; \; \ds\int_{0}^{\infty} N(x)\ud x=1, \quad \lambda_{0}>0,  \\
  \vspace{0.01cm}
   cg(x)\ds\frac{\partial}{\partial_{x}}\phi(x)-(B(x)+\lambda_{0})\phi(x)=-2B(x)\ds\int^{x}_{0}\kappa(y,x)\phi(y)\ud y, \\
   \vspace{0.01cm}
   \phi(x)\geq0; \; \ds\int_{0}^{\infty} \phi(x)N(x)\ud x=1. 
   \end{array} \right.\vspace{0.2cm}
\end{equation}
The use of this new problem allows to restrict the need for information to a non-temporal measure, and the problem becomes: How to recover information on $g,$ $B$ and $\kappa$ from an experimental measurement of the asymptotic profile $N$ and the global exponential rate of growth $\lambda_0$ of the population \footnote{Growth can naturally be balanced by death, by the addition for instance of a death term $d(x)N(x)$ on the left-hand side of the equation. This would lead to possible nonpositive rates $\lambda_0,$ but our whole study would remain unchanged.} ?

In the case when the equation models cell-division, direct measures of the growth rate $g(x)$ is possible. Direct measures of $\kappa$ is also possible, by a study of the sizes of the two daughter cells born from a mother. The most delicate point is thus the measure of the division rate $B,$ what implies to follow each cell from its birth to its division or death.

In \cite{P.Z} and \cite{D.P.Z}, the problem of recovering the division rate from a measured $N$ was addressed in the case when the growth rate is constant, \emph{i.e} $g(x)\equiv 1,$ and the daughter cells are twice smaller than their mother, \emph{i.e.} when  $\kappa(x,y)=\delta_{x=y/2}.$ In this case, Problem \eqref{(P)} writes:
\begin{equation}
c \ds\frac{\partial}{\partial_x} N + \big(B(x)+\lambda_0\big) N = 4 B(2x) N(2x).
\end{equation}
In this particular case, the inverse problem reads: How to recover $H=BN(x)$ solution of 
\begin{equation}\label{eq:inv:mitose}
{\cal L} (H)=F(N),\end{equation}
with ${\cal L}: \; H\to 4 H(2x)-H(x),$ and $F(N):= c\partial_x N + \lambda_0 N$ ?
 The method used to solve Equation \eqref{eq:inv:mitose} \footnote{the method was first developed in \cite{P.Z}, then investigated deeper and solved numerically in \cite{D.P.Z} in a deterministic setting, and in \cite{DHRR} in a statistical setting. It was also successfully applied to experimental data in \cite{D.M.Z}.}
 strongly uses the analytical study of the operator $\cal L,$ and it was shown that the most efficient technique was then to view the problem as written in the variable $y=2x$ rather than in $x$ (see the discussion in \cite{D.P.Z}). 

In this paper,  we address the inverse problem of determining the cell division rate $B$ when $g$ and $\kappa$ are known - or guessed - functions, but fully general ; hence, we cannot apply anymore the inversion of the operator $\cal L$ as done in \cite{D.P.Z}, and new tools have to be designed.

We model the experimental measure of the distribution $N$ by an approximation data $N_{\varepsilon}$ of $N$ satisfying $\Vert N-N_{\varepsilon}\Vert\leq\varepsilon$ for a suitable norm $\Vert $\textbullet $\Vert$ \footnote{A more precise model for the measured data, in a statistical setting, can be found in \cite{DHRR}.} 

The paper is organized as follows. We first study the regularity of the direct problem, what is a necessary step for a better understanding of the inverse problem. In a second part, we investigate the inverse problem of determining $B$ by the Quasi-reversibility and Filtering methods proposed in \cite{P.Z} and \cite{D.P.Z} and properly adapted to our general context. In a third part we develop new numerical approaches in order to recover the rate $B$ following the two regularization methods ; we give some numerical illustrations of our methods.
 
\subsection*{Main notations and assumptions}
We use the following notations.
\begin{equation}\label{def:P}
\mathbb{P}:=\bigl\{f\geq0\,:\,\exists\,\mu,\nu\geq0,\ \limsup_{x\to\infty}x^{-\mu}f(x)<\infty,\; \liminf_{x\to\infty}x^\nu f(x)>0\bigr\},\end{equation}
\begin{equation}\label{def:Lp0}
 L^p_0:=\bigr\{f,\ \exists a>0,\ f\in L^p(0,a)\bigl\},\qquad L^2_p:=L^2(\R_+,x^pdx).
\end{equation}

We work under the following technical assumptions, that guarantee well-posedness of Problem \eqref{(P)} as stated in \cite{D.G} (we refer to that paper for a complete discussion and justification).  
\begin{equation}\label{assum0}
\exists \, \,0<c < \frac{1}{2},\qquad \forall\; p\geq 2,\qquad D_p:=\int\limits_0^y \frac{x^p}{y^p} \kappa(x,y) dx \leq c <\frac{1}{2}.
\end{equation}
\begin{equation} \label{assum1}
B\in L^{1}_{loc}(\mathbb{R}^{*}_{+})\cap\mathbb{P}, \quad \exists \; \alpha_{0}\geq0, \; g\in L^{\infty}_{loc}(\mathbb{R_{+}},x^{\alpha_{0}}dx) \cap \mathbb{P} 
\end{equation}
\begin{equation}\label{assum2}
\forall \,K \; \text{compact in} \; ]0,+\infty[, \,\exists \; m_k>0 : g(x)\geq m_k \; \forall x\in K
 \end{equation}

\begin{equation}\label{assum3} 
\exists \, b\geq0, \;  supp B=[b,+\infty)
 \end{equation}

\begin{equation}\label{assum4} 
\exists \; C>0,\;  \gamma\geq 0, \frac{x^{\gamma}}{g(x)}\in L^1_0\; ; \; \int^{x}_{0}\kappa(z,y)\ud z \leq \min\big(1, C(\frac{x}{y})^{\gamma}\big) \; 
 \end{equation} 
\begin{equation}\label{assum5}
\frac{B(x)}{g(x)} \in L^{1}_{0} \quad ; \; \frac{xB(x)}{g(x)}\xrightarrow[x\to+\infty]{} +\infty.
 \end{equation}

\section{Regularity of the direct problem}
Before studying the inverse problem, it is necessary to have a proper knowledge of the direct one, which states as follows: What is the regularity  of the map $\Gamma: \; (c,B) \; \to \; (\lambda_0,N)$ solutions of Problem \eqref{(P)} ? How can we define a proper definition domain for $\Gamma$ ?

In \cite{D.P.Z}, Theorems 3.1. and 3.2 establish that the map $\Gamma_0: B\to (\lambda_0,N)$ is Lipschitz-continuous for $c=1$  fixed, $g=1,$ $\kappa=\frac{1}{2} \delta_{x=\f{y}{2}}$ and division rates $B$ such that $0<B_m \leq B\leq B_M <\infty;$ in other words, for division rates uniformly positive and uniformly bounded. 

In this paper, we want to state such results for general growth rates $g$ and division kernels $\kappa,$ with division rates $B$ not necessarily uniformly bounded.  Our study is thus first based on the well-posedness of this general eigenvalue problem \eqref{(P)}, as performed in \cite{D.G}. 

Let us first settle a proper definition space for the division rates $B.$ Theorem $1$ of \cite{D.G} states that, under Assumptions \eqref{as:kappa} and \eqref{assum0}-\eqref{assum5}, there exists a unique eigenpair $(\lambda_0, N, \phi)$ solution of Problem \eqref{(P)}. Hence, we first need that $g$ and $\kappa$ satisfy Assumptions \eqref{as:kappa}, \eqref{assum0}--\eqref{assum2}. Then, to study the regularity of the map $\Gamma: (c,B)\to (\lambda,N),$ one needs not only that such division rates $B$ satisfy Assumptions \eqref{assum1}, \eqref{assum3} and \eqref{assum5} but also that they satisfy them \emph{uniformly}. This leads to the following definition.
\begin{Def} Let $g,\,\kappa$ satisfying Assumptions \eqref{as:kappa}, \eqref{assum0}-\eqref{assum2}. For a constant $b\geq 0$ and functions $f_0\in L_0^1,$ $f_\infty \xrightarrow[x\to+\infty]{} \infty,$ one defines the set
$${\cal D}(b,f_0,f_\infty) :=\left\{ B \in L^\infty_{loc} (\R_+^*) \cap \mathbb{P}, \quad Supp(B)=[\tilde b \leq b, +\infty), \quad \frac{B}{g} \leq f_0,\quad \frac{x B}{g} \geq f_\infty
\right\}.$$
\end{Def}
In such a set, division rates $B$ satisfy uniformly Assumption \eqref{assum5}, what allows to use the powerful estimates proved in \cite{D.G}. 

Under such assumptions, we also recall that we have the following results (see Theorem 1 in \cite{D.G}) for the unique solution $(\lambda_0,N,\phi)$ to Problem \eqref{(P)}:  
  \begin{equation}x^\alpha gN\in L^p(\mathbb{R+}) \;\;  \forall \;  \alpha\geq -\gamma, \; \forall 1\leq p \leq +\infty \; ; \quad x^\alpha gN \in W^{1,1}(\mathbb{R+}) \; \;  \forall \alpha\geq0, \label{prop:N}\end{equation}
  and 
  \begin{equation}\label{prop:phi} \exists \; k>0, \;C>0,\; \phi(x)\leq C(1+x^k) ; \quad g\partial_x \phi \in L^{\infty}(\mathbb{R+}).\end{equation}
The two following fundamental estimates are straightfully obtained by integration on $[0,\infty[$ of Equation \eqref{(P)} or \eqref{(P)} multiplied by  $x$:
 \begin{equation}\label{estim1}
 \lambda_0=\int^{+\infty}_{0} B(x)N(x)\ud x,
 \end{equation}    
  
 \begin{equation}\label{estim2}
 \int^{+\infty}_{0}xN(x)\ud x=\frac{c}{\lambda_0}\int^{+\infty}_{0} g(x)N(x)\ud x.
 \end{equation}
 
We are now ready to state our regularity result.

\begin{theor}\label{cont_th}
Let parameters $g$ and $\kappa$  satisfy Assumptions \eqref{assum0}-\eqref{assum2}, then \\
\begin{itemize}
\item[i)] The map $\Gamma: (c,B)\longmapsto (\lambda_{0},N)$
is:

\textbullet \;  continuous in $(c,B)$ under the $L^\infty-$ weak-\textasteriskcentered topology for $B\,$
 from any set $\mathbb{R_{+}^*}\times {\cal D}(b,f_0,f_\infty)$  to  $\mathbb{R_{+}^*}\times L^{1}\cap L^{\infty}(\mathbb{R_{+}}) $. \\
 \\
 \textbullet \; injective.\\ 
 \item[ii)] Let moreover $g$ satisfy $\f{x^\gamma}{g} \in L^2_0$ with $\gamma$ defined in Assumption \eqref{assum4}. Then
 the map $\Gamma$ is Lipschitz-continuous under the strong topology of $\R_+^* \times L^2 \cap {\cal D}(b,f_0,f_\infty).$ More precisely,  denoting $\delta B=\bar{B}-B, \; \Delta= \Vert \bar{B}-B\Vert_{L^{2}(\mathbb{R_{+}})}$, $\delta c = \bar{c} -c,$ $\delta N=\bar{N}-N,$ $\delta \lambda=\bar{\lambda_{0}}-\lambda_{0},$ we have the following estimates, for $k$ as in \eqref{prop:phi}: 
 $$\vert \delta \lambda\vert\leq C_1(B,\bar{B})\Delta + C_2(B,\bar B )|\delta c|,\quad \Vert \delta N\Vert_{L^2 (\R_+)} \leq C_3 (B,\bar{B}) \Delta,\quad  \;$$
 with 
 $$  C_1=C \ds \Vert \frac{\phi}{1+x^k}   \Vert_{L^\infty} \frac{\Vert (1+x^{k})\bar{N}\Vert_{L^{2}}}{\int^{\infty}_{0}\bar{N}\phi\ud x},\quad C_2=\frac{\lambda_0+C}{c}\Vert (1+x^{k})g\bar{N}\Vert_{L^{1}(\mathbb{R_{+}})} + 
 \Vert (1+x^{k})g\bar{N}\Vert_{L^{2}} \Vert B\Vert_{L^{2}},$$ 
where $C>0$ is an absolute constant.

\end{itemize}  
\end{theor} 

\begin{proof}\\
\begin{itemize}
\item[i)]  The continuity of the map $\Gamma$ directly follows from the proof given  in \cite{D.G}, Theorem 1. Therefore, we only sketch the main steps and let the reader refer to this article. 

 Let $c_n \to c >0 $ in $\R_+^*$ and $B_n \overset{*}{\rightharpoonup} B$ in ${L^{\infty}(\R_+)}.$  Denoting $(\lambda_n,N_n)$ the respective eigenpairs solutions of Problem \eqref{(P)} settled for $(c_n,B_n),$  we can prove the same uniform estimates for $N_n$ as in \cite{D.G} due to the fact that since $B_n \in {\cal D}(b,f_0,f_\infty),$ Assumption \eqref{assum5} is uniformly verified. Such estimates give strong compactness in $L^1$ for $N_n,$ and hence, up to a subsequence, we have a strong convergence of $(\lambda_n,N_n)$ to $(\lambda, N).$ Similarly, we prove $\lambda >0,$ and passing to the limit in the equations for $N_n,$ we deduce that $(\lambda, N)$ has to be the solution of Problem \eqref{(P)} settled for $(c,B).$ Since such a solution is unique,  the whole sequence $(\lambda_n, N_n) $ converges to it. 

Let us show by contradiction that $\Gamma$ is an injection function.\\
 Let $B_{i}\in L^{1}_{loc}(\mathbb{R^{*}_{+}})$ and $c_{i}$ positive constants $\forall i\in \{1,2\}$ such that $(c_{1},B_{1})\neq (c_{2},B_{2})$ and $\Gamma(c_{1},B_1)=\Gamma(c_{2},B_{2})=(\lambda_0,N)$. \\ 
We then integrate the two equations satisfied by $(\lambda_0,N)$ against the weight $x,$ to obtain
 $$\int^{\infty}_{0}g(x)N(x)\ud x = \frac{\lambda_{0}}{c_{1}}\int^{\infty}_{0}yN(y)\ud y = \frac{\lambda_{0}}{c_{2}}\int^{\infty}_{0}yN(y)\ud y, $$
what implies $c_1=c_2.$
By the contradiction assumption we get $B_{1}\neq B_{2}$, so by making the difference between the following Equations \eqref{eq3}, \eqref{eq4} with consideration to the equality $c_{1}=c_{2}$
\begin{eqnarray}
c_{1}\f{\p}{\partial_{x}}(g(x)N(x))+(B_{1}(x)+\lambda_0)N(x)&=&2\int_{0}^{\infty}B_{1}(y)\kappa(x,y)N(y)\ud y, \label{eq3}\\
c_{2}\f{\p}{\partial_{x}}(g(x)N(x))+(B_{2}(x)+\lambda_0)N(x)&=&2\int_{0}^{\infty}B_{2}(y)\kappa(x,y)N(y)\ud y, \label{eq4}
\end{eqnarray}
  we obtain after multiplying by $x^{p}$, $p\geq 2$ the following relation
 $$ x^{p}\vert B_{1}-B_{2}\vert(x)N(x)\leq 2\int_{0}^{\infty}x^{p}\vert B_{1}-B_{2}\vert(y)N(y)\kappa(x,y)\ud y. $$
 We integrate this relation on $(0,\infty)$ and due to Assumption \eqref{assum0} for $p\geq2$ we deduce the following strict inequality :
 $$\int_{0}^{\infty} x^{p}\vert B_{1}-B_{2}\vert(x)N(x)\ud x < \int_{0}^{\infty}x^{p}\vert B_{1}-B_{2}\vert(x)N(x)\ud x ,\quad \forall p\geq2$$
what is contradictory. \\
 \item[ii)] First, the fact that $\ds\f{x^\gamma}{g} \in L^2_0$ implies that for all $p\geq 0,$ $N\in L^2\big((1+x^p)\ud x\big).$ Indeed, by \eqref{prop:N} and since $g \in \mathbb{P},$ $\ds\int^\infty N^2 (x) (1+x^p)^2 \ud x <\infty,$ and it only remains to bound $\ds\int_0 N^2 (1+x^p)^2 \ud x.$ This is given by writing
$N^2 (x) = (N^2 {g^2} x^{-2\gamma})( \ds\f{x^{2\gamma}}{g^2}),$ product of an $L^\infty$ function with a $L^1_0$ function.

 By making the sum between the two following equations
 \begin{equation*}
\begin{array}{ll}
{\color{red}\phi(x)}\bar c\ds\f{\p}{\partial_{x}}(g(x)\bar{N}(x) ) +{\color{red}\phi(x)}(\bar{B}(x)+\bar{\lambda}_{0})\bar{N}(x)={\color{red}\phi(x)}2\ds\int^{\infty}_{0}\bar{B}(y)\kappa(x,y)\bar{N}(y) \ud y \\
\\
{\color{red}\bar{N}}(x)cg(x)\ds\f{\p}{\partial_{x}}\phi(x) -{\color{red}\bar{N}(x)}(B(x)+\lambda_{0})\phi(x)=-2{\color{red}\bar{N}(x)}B(x)\ds\int^{x}_{0}\kappa(y,x)\phi(y) \ud y 
 \end{array} \vspace{0.05cm}
\end{equation*}
we obtain
\begin{eqnarray*}
\delta c \; \phi \ds\f{\p}{\partial_x}(g\bar N)  + \f{\p}{\partial_x}\bigl(cg\bar{N}\phi\bigr)(x) +\bigl(\phi\bar{N}[\bar{B}+\bar{\lambda}_{0}-B-\lambda_0]\bigr)(x)=2\phi(x)\int^{\infty}_{0}\bar{B}(y)\kappa(x,y)\bar{N}(y) \ud y \\
-2\bar{N}(x)B(x)\int^{x}_{0}\kappa(y,x)\phi(y) \ud y
\end{eqnarray*}
we then integrate this equation on $[0,\infty)$ that leads
$$\delta c \; \phi \f{\p}{\partial_x}(g\bar N) dx +\delta\lambda\int^{\infty}_{0}\phi\bar{N}\ud x +\int^{\infty}_{0}\phi\bar{N}\delta B\ud x=2\int^{\infty}_{0}\delta B(y)\bar{N}(y)\biggl(\int^{\infty}_{0}\phi(x)\kappa(x,y)\ud x\biggr)\ud y .$$  
So $$\delta\lambda\int^{\infty}_{0}\phi(x)\bar{N}(x)\ud x =\int^{\infty}_{0}\delta B(x)\bar{N}(x)\biggl(2\int^{\infty}_{0}\phi(y)\kappa(y,x)\ud y -\phi(x) \biggr)\ud x + \delta c \int g\bar N \f{\p}{\partial_x}\phi \ud x  .$$
The first term of the left-hand side gives the term with $C_1(B,\bar B)$ of the estimate on $\delta \lambda$ by using the fact that $ \exists \; k>0, \ds\frac{\phi}{1+x^k} \in L^\infty(\mathbb{R}_{+})$. For the second term, we use the equation for $\phi$ and write
$$ \delta c \int g\bar N \f{\p}{\partial_x}\phi \ud x = \f{\delta c}{c} \int g\bar N \biggl((B+\lambda_{0})\phi-2B\ds\int^{x}_{0}\kappa(y,x)\phi(y)\ud y \biggr)dx,
$$  
and it provides the term with $C_2(B,\bar B)$ in the estimate for $\delta \lambda.$
\newline

To prove the estimate on $\delta N,$ we make the difference between the two following equations
\begin{equation*}
\begin{array}{ll}
\vspace{0.2cm}
c\ds\frac{\partial}{\partial_{x}}(g(x)\bar{N}(x) ) +(\bar{B}(x)+\bar{\lambda}_{0})\bar{N}(x)=2\ds\int^{\infty}_{0}\bar{B}(y)\kappa(x,y)\bar{N}(y) \ud y \\
c\ds\frac{\partial}{\partial_{x}}(g(x)N(x) ) +(B(x)+\lambda_{0})N(x)=2\ds\int^{\infty}_{0}B(y)\kappa(x,y)N(y)\ud y \\
 \end{array} \vspace{0.05cm}
\end{equation*}
we obtain
\begin{eqnarray*}
\delta c \frac{\partial}{\partial_x} (g\bar N) + c\ds\f{\p}{\partial_{x}}(g\delta N)&+&\biggl((\bar{\lambda}_{0}+\bar{B})\bar{N}-(\lambda_{0}+B)N\pm(\lambda_{0}+B)\bar{N}\biggr)
\\&=&2\int^{\infty}_{0}\biggl(\bar{B}\bar{N}-BN \pm B\bar{N}\biggr)(y)\kappa(x,y)\ud y .
\end{eqnarray*}
That implies
\begin{eqnarray*}
\delta c \ds\frac{\partial}{\partial_x} (g\bar N) +c\ds\frac{\partial}{\partial_{x}}(g\delta N)+(\lambda_{0}+B)\delta N&=&\biggl[2\int^{\infty}_{0}\bar{N}(y)\kappa(x,y)\delta B(y) \ud y -(\delta\lambda+\delta B)\bar{N} \biggr] \\
&&+2\int^{\infty}_{0}B(y)\delta N(y)\kappa(x,y)\ud y. 
\end{eqnarray*}
We recast the previous equation as follows
\begin{equation}\label{eq1}
c\ds\frac{\partial}{\partial_{x}}(g(x)\delta N(x))+(\lambda_{0}+B(x))\delta N(x)=2\int^{\infty}_{0}B(y)\delta N(y)\kappa(x,y)\ud y +\delta R(x),
\end{equation}
with 
\begin{equation}\label{def:deltaR}
\delta R(x)= 2\int^{\infty}_{0}\bar{N}(y)\kappa(x,y)\delta B(y) \ud y -(\delta\lambda+\delta B)\bar{N} - \delta c \f{\p}{\partial_x} (g\bar N) 
\end{equation}
We can bound $\Vert \delta R(x)\Vert _{L^2}$ as we previously bound $|\delta \lambda|.$ The estimate on $\Vert \delta N \Vert_{L^2}$ thus follows from the following lemma.
 \end{itemize}\end{proof}

\begin{lem}
Under the assumptions of Theorem \ref{cont_th} for $g$ and $\kappa,$ with $\delta N$ defined as in  Theorem \ref{cont_th} and $\delta R$ defined by \eqref{def:deltaR}, there exists $\nu(c,B) >0$ a constant depending only on the eigenvalue problem \eqref{(P)} stated for given parameters $c >0$ and $B\in L^2\cap {\cal D}(b,f_0,f_\infty)$ such that, for all $\bar c \geq c_0>0$ and $\bar B \in L^2\cap {\cal D}(b,f_0,f_\infty),$ one has
$$  \nu\Vert \delta N\Vert_{L^2(\mathbb{R}_{+})}\leq \Vert\delta R\Vert_{L^2(\mathbb{R}_{+})}.$$
 \end{lem}
 \begin{proof}\\
We argue by contradiction and assume that for a sequence $c_k \geq c_0>0,$ $\bar{B}_{k} \in L^2\cap {\mathcal D} (b,f_0,f_\infty)$, one has, for a vanishing sequence $\nu_k,$
$$  \nu_k \Vert \delta N_k\Vert_{L^2(\mathbb{R}_{+})}\geq \Vert\delta R_k\Vert_{L^2(\mathbb{R}_{+})},$$
with $\delta N_k=\bar N_k -N,$ $\bar N_k$ solution of Problem \eqref{(P)} stated for $\bar c_k$ and $\bar B_k,$ $\delta R_k$ defined by \eqref{def:deltaR} stated for $\delta N_k.$ 

As for the proof of continuity above, compactness arguments as done in \cite{D.G} lead us to extract a converging subsequence $\bar N_k \to \bar N$ strongly in $L^1,$ so $\delta N_k \to \delta N$   strongly in $L^1.$ Moreover, estimates as in \cite{D.G} imply that $\bar N_k$ is uniformly bounded in $L^2$ (we write $\ds\bar N_k^2= x^{-2\gamma} g^2 N^2\f{x^{2\gamma}}{g^2}$ and use the assumption $\ds\f{x^\gamma}{g} \in L^2$  together with the result \eqref{prop:N}, result which is uniform for all $\bar N_k$), hence $\delta N_k$ satisfy Equation \eqref{eq1} with $\Vert \delta R_k \Vert_{L^2} \to 0.$ Passing to the limit, it implies that $\delta N$ satisfies Equation \eqref{(P)}, so by uniqueness of a solution we have $\delta N = C N$ for a given constant $C\in \R.$ Since $\int N dx=\int \bar N_k dx=1,$ we have $\int \delta N dx = C =0:$ it is contradictory with our assumption on $(\nu_k).$\qed

 \end{proof}

\section{The inverse problem and its regularization}

As in \cite{D.M.Z,DHRR}, we consider the problem of recovering the cell division rate $B$ and the constant $c$ from the \emph{a priori} knowledge of the shape of the growth rate $g(x)$ and the experimental measure of the asymptotic distribution $N$ and exponential growth $\lambda_0.$ To model this, we suppose that we have  two given measurements $N_\eps \in L^1 \cap L^\infty(\R_+)$ and $\lambda_\eps>0$ such that 
$||N-N_\eps||_{L^2\big((1+x^p)\ud x\big)} \leq \eps,$  $|\lambda_0 - \lambda_\eps|\leq \eps.$\footnote{See \cite{DHRR} for a statistical viewpoint on the data $(N_\eps, \lambda_\eps):$ supposing that $N_\eps \in L^2$ means that we deal with some preprocessed data. However, once the problem is solved in a deterministic setting, as we do in this article, it is immediate to apply the method of \cite{DHRR} to this general case.} The problem is: How to get estimates $(c_\eps,B_\eps)$ of $(c,B)$ solutions of
\begin{equation} \label{statio}
c\ds\frac{\partial}{\partial_{x}}(g(x)N(x) ) +(B(x)+\lambda_{0})N(x)=2\int^{\infty}_{0}\kappa(x,y)B(y)N(y) \ud y.
\end{equation} 
First, it is clear that $B$ cannot be recovered from Equation \eqref{statio} when the distribution $N$ vanishes: our inverse problem consists in recovering $H=BN$ rather than $B$ directly.  Our problem can now be viewed as: How to recover $(c,H)$ solution of
\begin{equation} \label{statio:inverse}
{\cal L}_\kappa ( H) (x):=H(x) - 2\int^{\infty}_{0}\kappa(x,y)H(y) \ud y
= - c\ds\frac{\partial}{\partial_{x}}(g(x)N(x) ) -\lambda_{0} N(x)
\end{equation} 
when we have measurements $(\lambda_\eps, N_\eps)$ of $(\lambda_0,N)$?

Secondly, since the measure $N_\eps$ is supposed to be in $L^2,$ there is no way of directly controlling $\ds\frac{\partial}{\partial_{x}}(gN_\eps)$ even if $g$ is known (see Section 2 of \cite{P.Z} for a discussion, or yet \cite{Engl}).
 To overcome this difficulty, two regularization methods were proposed in \cite{P.Z, D.P.Z} for the particular case of division into two equal cells, \emph{i.e.} when $\kappa(x,y)=\delta_{x=y/2},$ a third method has also been proposed in \cite{Groh}, and a statistical treatment to estimate the derivative in \cite{DHRR}. Indeed, looking at the problem in terms of $H=BN$ and not in terms of $B$ makes it \emph{almost} linear in $H;$ almost, because $\lambda_0$ being also measured, the term $\lambda_0 N$ can be viewed as quadratic. Hence, the classical tools designed to regularize linear inverse problems (see \cite{Engl}) can be used, as illustrated by the  three foreseen methods, as soon as the operator ${\cal L}_\kappa$ can be inverted.

This is the third and last difficulty: inverse the operator ${\cal L}_\kappa$ defined by Equation \eqref{statio:inverse}. None of the three regularization methods of \cite{P.Z, D.P.Z, Groh} can be directly applied here: indeed, they strongly used the fact that for the kernel $\kappa=\delta_{x=\f{y}{2}},$ the left-hand side of Equation \eqref{statio:inverse} simplifies in $4BN(2x)-B(x),$ and can be viewed as an equation written in $y=2x.$ Then, a central point of the proofs in \cite{P.Z} as well as in \cite{D.P.Z} or \cite{Groh} is the use of the Lax-Milgram theorem for the coercitive operator ${\cal L}: H\to 4H(y) - H(\f{y}{2}).$ 

Nothing such as that can be written here, and the main difficulty, numerically as well as theoretically, is to deal with a nonlocal kernel $\int \kappa(x,y) H(y)dy.$  The operator $\cal L$ is replaced by ${\cal L}_\kappa.$
For $\kappa=\delta_{x=\f{y}{2}},$ ${\cal L}_\kappa$ has been proved in \cite{D.P.Z} (Proposition A.1. in the appendix) to be coercitive in $L^2(x^p dx)$ if $p>3,$ or in contrary ${\cal L}$ is coercitive if $p<3.$ Due to the nonlocal character of the kernel, it seems more natural now to look for cases when the first part of the operator ${\cal L}_\kappa,$ \emph{i.e.} identity, dominates the nonlocal part $2\int\limits_x^\infty u(y) \kappa(x,y) \ud y.$ This is expressed by the following proposition.

\begin{prop}\label{prop:filt1}
Let $\kappa$ satisfy Assumption \eqref{as:kappa}. For $r,q\geq 0,$ we define the following quantities:
\begin{equation}\label{def:CkDq}
C_r:=\sup_x \; \int\limits_x^\infty  \f{x^r}{y^r} \kappa(x,y) dy,\qquad D_q:=\sup_y \; \int\limits_0^y \f{x^{q}}{y^{q}} \kappa(x,y) dx.
\end{equation}
If $0\leq r\leq p$ are such that
\begin{equation}\label{coer_cond}
C_r D_{p-r} <\frac{1}{4},
 \end{equation}
Then for all $F \in L^2(\R_+, x^p \ud x)$ there exists a unique solution $u \in L^2(x^p\ud x)$ to the following problem:
\begin{equation}
u(x) - 2\ds\int\limits_x^\infty u(y) \kappa(x,y) \ud y = F,
\end{equation}
and we have the following estimate
$$||u||_{L^2(x^p \ud x)} \leq \frac{1}{1-2\sqrt{C_r D_{p-r}}}||F||_{L^2(x^p \ud x)}.$$
\end{prop}
\begin{proof}
We define the bilinear form
\begin{eqnarray*}
\mathcal{A}(u,v)=\int^{+\infty}_{0}u(x)v(x)x^p \ud x-2\int^{+\infty}_{0}u(x)\int^{+\infty}_{x}\kappa(x,y)v(y)\ud y\, x^p \ud x=<u,v>_{L^2(x^p \ud x)} - {\mathcal{B}}(u,v), 
\end{eqnarray*}
where $<,>_{L^2}$ denotes the scalar product. We apply the Lax-Milgram theorem in $L^2(x^p \ud x).$ Indeed, we have by Cauchy-Schwartz, for any constant $C>0:$
$$\begin{array}{lll}
{\mathcal{B}} (u,v) &=&2\ds\int_0^\infty \int_x^\infty {x^p} u(x)  v(y) \kappa(x,y) \ud y \ud x 
=2\ds\int_0^\infty \int_0^y {x^p} u(x)  v(y) \kappa(x,y) \ud x \ud y
\\ \\
&\leq& \ds\int\limits_0^\infty\int\limits_0^y \big(C u^2(x) x^p \f{x^r}{y^r} + \f{1}{C} u^2(y) y^p \f{x^{p-r}}{y^{p-r}}\big) \kappa(x,y) dx dy \\ \\ 
&=&C \ds\int\limits_0^\infty u^2(x) x^p \big(\ds\int\limits_x^\infty  \f{x^r}{y^r} \kappa(x,y) dy \big) dx + \f{1}{C} \ds\int\limits_0^\infty u^2(y) y^p \big(\ds\int\limits_0^y \f{x^{p-r}}{y^{p-r}} \kappa(x,y) dx \big) dy \\
&\leq& (C C_r + \f{D_{p-r}}{C}) \ds\int\limits_0^\infty u^2(x) x^p dx,
\end{array}
$$ 
The minimum is reached for $C=\sqrt{\f{C_k}{D_{p-k}}},$ So finally we have
$${\mathcal{B}} (u,v) \leq 2 \sqrt{C_r D_{p-r}} ||u||_p^2,$$
what proves the continuity of the bilinear forms $\mathcal B$ and $\mathcal A.$ Moreover,  it implies
$$\mathcal{A} (u,u) \geq (1 - 2 \sqrt{C_r D_{p-r}}) ||u||^2_{L^2(x^p \ud x)}=\beta ||u||^2_{L^2(x^p \ud x)},$$
with $\beta=1-2\sqrt{C_r D_{p-r}} >0$ under assumption \eqref{coer_cond}. It ends the proof of Proposition \ref{prop:filt1}.
\end{proof}
\begin{Rmq}
Assumption \eqref{coer_cond} 
can be linked to Assumption \eqref{assum0}. For self-similar kernels
$\kappa(x,y)=\ds\f{1}{y} k_0 (\f{x}{y})$ with $\ds\int_0^1 k_0(z)dz=1,$ defining $I_r=\int\limits_0^1 z^r k_0(z) dz$ we obtain
$C_r=I_{r-1}$ and $D_q=I_q,$ so that Assumption \eqref{coer_cond} is reduced to 
$$I_{r-1}I_{p-r} <\f{1}{4}.$$
For the equal mitosis kernel $\kappa(x,y)=\delta_{x=\f{y}{2}},$ since $I_q=2^q,$ Assumption \eqref{coer_cond} is verified for $p>3:$ we recover part of the result of the proposition of \cite{D.P.Z}. It corresponds to the cases when the first part of the bilinear form (\emph{i.e.}, $\ds\int uvx^p \ud x$) dominates the second one ($\ds\iint\kappa(x,y) u(x)v(y) x^p \ud x \ud y$).

For the uniform kernel $k_0(z)=1,$ the equality $I_q=\f{1}{q+1}$ here again implies that Assumption \eqref{coer_cond} is verified for $p>3.$ 

More generally, for these self-similar kernels Assumption \eqref{assum0} is equivalent to  $I_{r\geq 2}<\f{1}{2}.$ Since $I_0=1$ and $I_1=\f{1}{2}$ by Assumption \eqref{as:kappa}, taking $r=2$ in Assumption \eqref{coer_cond} leads to
$$I_{p-2} <\f{1}{2},$$
and due to Assumption \eqref{assum0} this is always true for $p\geq 4.$
\end{Rmq}

\subsection{Filtering method}

This regularization method consists in looking for a solution $H_{\eps,\alpha}$ of the following regularized problem
\begin{equation} \label{Pestf}
{\cal L}_\kappa(H_{\varepsilon,\alpha})(x):=H_{\varepsilon,\alpha} (x) - 2\ds\int^{+\infty}_{0}\kappa(x,y)H_{\varepsilon,\alpha}(y)\ud y 
= \rho_\alpha * \biggl(- c_{\varepsilon,\alpha}\frac{\partial}{\partial_{x}}\bigl(g(x)N_{\varepsilon}(x)\bigr) -\lambda_{\varepsilon}N_{\varepsilon}(x)\biggr),
 \vspace{0.2cm}\\
\end{equation}
where $\rho_{\alpha}$ is a mollifiers sequence defined by 
\begin{equation}\label{def:mollifier}
\rho_{\alpha}(x)=\ds\frac{1}{\alpha}\rho(\frac{x}{\alpha}), \quad \rho\in\mathbb{C}^{\infty}_{c}(\mathbb{R}), \quad \int_{0}^{\infty}\rho(x)\ud x=1, \quad \rho\geq0, \quad Supp(\rho)\subset [0,1].
\end{equation} \\
One notices that $c_{\eps,\alpha}$ is uniquely defined: indeed, integrating Equation \eqref{Pestf} against the weight $x$ leads to 
\begin{equation} \label{def:cepsalpha}
c_{\eps,\alpha} = \lambda_\eps \ds\f{\ds\int x N_\eps \ud x}{\ds\int \rho_\alpha * (g N_\eps) \ud x} .  \end{equation}
We want to study the well-posedness of this problem and estimate the distance between $B_{\eps,\alpha}=\frac{H_{\eps,\alpha}}{N_{\eps,\alpha}}$ and $B$ in order to choose an optimal approximation rate $\alpha.$ This is given by the following result.

\begin{theor}\label{th:filter}
Let $g,$ $B$ and $\kappa$ satisfy Assumptions \eqref{as:kappa} and \eqref{assum0}--\eqref{assum5}, and moreover $\ds\f{x^\gamma}{g} \in L^2_0$ with $\gamma$ defined in Assumption \eqref{assum4}. Let $(\lambda_0,N)$ the unique eigenpair solution of Problem \eqref{(P)} (as stated in \cite{D.G}).
Let $p>1$ satisfy Assumption \eqref{coer_cond}.
 Let $N_\eps \in L^1\cap L^\infty (\R_+)$ and $\lambda_\eps >0$ satisfy $\Vert g(N-N_\eps)\Vert_{L^2(x^p \ud x)} \leq \eps \Vert gN\Vert_{L^2(x^p \ud x)},$ $|\lambda_\eps -\lambda_0|\leq \eps \lambda_0,$ $\Vert N - N_\eps \Vert_{L^1\big((1+x+g(x))\ud x\big)} \leq \eps \Vert N\Vert_{L^1\big((1+x+g(x))\ud x\big)} $ and \\ $\Vert N_\eps - N\Vert_{L^2\big((x^p+1)\ud x\big)} \leq \eps \Vert N\Vert _{L^2\big((x^p+1)\ud x\big)}.$ 
 
Then there exists a unique solution $H_{\eps,\alpha} \in L^2(x^p \ud x)$ to Problem \eqref{Pestf}.

Defining $B_{\eps,\alpha}:=\chi_{N_{\eps,\alpha}(x) \neq 0} H_{\eps,\alpha} / N_{\eps,\alpha}$ we have the following estimates: 
\begin{equation}|c_{\eps,\alpha} - c| \leq C(p,\rho,N) ( \alpha + {\eps}),\label{est:filter:c}
\end{equation}
\begin{equation}||B_{\eps,\alpha} - B||_{L^2 (x^p N^2 \ud x)} \leq C(p,\rho,N) ( \alpha + \frac{\eps}{\alpha}),\label{est:filter}
\end{equation}
where $C$ is a constant depending on $p,$ moments of $\rho$ and $\f{\p}{\p_x} \rho,$ $\lambda_0,$ $\Vert gN\Vert_{H^2\big((1+x^p)\ud x\big)},$ $\Vert N\Vert_{L^1\big((1+x+g(x))\ud x\big)},$ $\Vert gN\Vert_{W^{1,1}(\ud x)}$ and $\Vert N\Vert _{H^1\big((x^p+1)\ud x\big)}.$
\end{theor}

The estimate \eqref{est:filter} of Theorem \ref{th:filter} relies, on the one hand, on the estimate of Proposition \ref{prop:filt1}, and, on the other hand, on general approximation properties of  the mollifiers, as expressed by Lemma \ref{lm:mollifier} right above. 
\begin{lem}\label{lm:mollifier}
Let $p>1,$ $f\in L^2 ((x^p +1)\ud x),$ $\rho_\alpha$ a mollifiers sequence defined by \eqref{def:mollifier} and $0<\alpha<1.$ Then we have the following estimates.
\begin{enumerate}
\item $||f * \rho_\alpha ||_{L^2(x^p \ud x)} \leq C(p,\rho) ||f||_{L^2\big((x^p+1) \ud x\big)},$ with $C(p,\rho)$ only depending on $p$ and moments of $\rho,$ \label{lm:mollifier:estim1}
\item \label{lm:mollifier:estim2}
$\ds||\f{\p }{\p x} (f * \rho_\alpha) ||_{L^2(x^p \ud x)} \leq \f{1}{\alpha}C(p,\rho)||f||_{L^2\big((x^p+1) \ud x\big)},$ with $C(p,\rho)$ only depending on $p$ and moments of $\rho$ and $\f{\p}{\p x} \rho.$
\item \label{lm:mollifier:estim3}
$||f*\rho_\alpha - f||_{L^2(x^p \ud x)} \leq C(\rho) \alpha ||f||_{H^1 (x^p \ud x)}$ if $f\in H^1 \big((1+x^p)\ud x\big)$ 
\item \label{lm:mollifier:estim4}
$||f*\rho_\alpha - f||_{L^1} \leq C(\rho)\alpha ||f||_{W^{1,1}}$
\item \label{lm:mollifier:estim5}
$||\rho_\alpha * f||_{L^1} \leq ||f||_{L^1}.$
\end{enumerate}
\end{lem}
\begin{proof}The proof of this result is classical and relies on Minkowski inequality for convolution products ; we let it to the reader. 
\end{proof}

{\bf Proof of Theorem \ref{th:filter}.}
We decompose the left-hand side of Estimate \eqref{est:filter} as follows
\begin{eqnarray*}
\Vert B_{\varepsilon,\alpha}N-BN \Vert_{L^{2}(x^{p}\ud x)}&=&\Vert B_{\varepsilon,\alpha}(N - N_\alpha + N_\alpha - N_{\varepsilon,\alpha})+ H_{\varepsilon,\alpha} - BN \Vert_{L^{2}(x^{p}\ud x)}\\
&\leq& \Vert B_{\varepsilon,\alpha} \Vert_{L^\infty} \biggl(\Vert N - N_{\alpha} \Vert_{ L^{2}(x^p \ud x)} +
\Vert N_\alpha - N_{\varepsilon,\alpha} \Vert_{ L^{2}(x^p \ud x)}\biggr)
+ \Vert H_{\varepsilon,\alpha} - BN \Vert_{L^{2}(x^p \ud x)} 
\end{eqnarray*}
On the right-hand side, the first term is bounded by $C(p,\rho)\alpha \Vert N\Vert_{H^1\big((x^p+1)\ud x\big)} $ due to Lemma \ref{lm:mollifier}, Estimate \ref{lm:mollifier:estim3}. The second term is bounded by $C(p,\rho) \eps \Vert N\Vert_{L^2\big((x^p+1)dx\big)}$ due to Lemma \ref{lm:mollifier}, Estimate \ref{lm:mollifier:estim1} applied to $f=N-N_\eps.$ For the third term, we apply Proposition \ref{prop:filt1} to $u=H_{\eps,\alpha} - BN$ and $F= \rho_\alpha * \biggl( c_{\varepsilon,\alpha}\ds\f{\p}{\partial_{x}}\bigl(g N_{\varepsilon}\bigr) +\lambda_{\varepsilon}N_{\varepsilon}\biggr) -  \biggl(c\ds\f{\p}{\partial_{x}}\bigl(g(x)N(x)\bigr) +\lambda_0 N\biggr) $. We treat these terms in a similar manner. Let us detail briefly the most binding term (with the notation $L^2_p=L^2(x^p \ud x)$):
$$\begin{array}{ll}\ds\Vert \rho_\alpha *  c_{\varepsilon,\alpha}\f{\p}{\partial_{x}}\bigl(g N_{\varepsilon}\bigr) - c\f{\p}{\partial_{x}}\bigl(gN\bigr)
\ds\Vert_{L^2_p} \leq &c_{\eps,\alpha} \biggl(\ds\Vert \f{\p}{\partial_x} \rho_\alpha * (gN_\eps - gN)\Vert_{L^2_p} +\ds\Vert \rho_\alpha*\f{\p}{\partial_{x}}\bigl(gN\bigr) - \f{\p}{\partial_{x}}\bigl(gN\bigr) \Vert_{L^2_p} \biggr)
\\
& + \ds |c_{\eps,\alpha}-c| \Vert \f{\p}{\partial_{x}}\bigl(gN\bigr)\Vert_{L^2_p}.\end{array}$$
The first term is bounded by $C \ds\frac{\eps}{\alpha} \Vert gN\Vert_{L^2(x^p \ud x)}$ by Lemma \ref{lm:mollifier} Estimate \ref{lm:mollifier:estim2}, the second one by $C \alpha \Vert gN\Vert_{H^2(x^p\ud x)}$ by Estimate \ref{lm:mollifier:estim3}. For the third term we write
$$\begin{array}{lll}
\vspace{0.15cm}
\vert c_{\eps,\alpha} -c\vert &= &\vert\lambda_\eps \ds\f{\int x N_\eps \ud x}{\int \rho_\alpha * g N_\eps \ud x} - \lambda_0 \f{\int x N(x) \ud x}{\int g(x)N(x) \ud x}\vert \\
\vspace{0.15cm}
&=&\ds|\lambda_\eps \f{\int x N_\eps}{\int \rho_\alpha * g N_\eps \ud x} \pm \lambda_\eps \f{\int x N\ud x}{\int \rho_\alpha * g N_\eps\ud x} \pm \lambda_0 \f{\int x N \ud x}{\int \rho_\alpha * g N_\eps \ud x} - \lambda_0 \f{\int x N \ud x}{\int g N \ud x} | \\
\vspace{0.15cm}
&\leq & \ds |\lambda_\eps \f{\int x N_\eps}{\int \rho_\alpha * g N_\eps \ud x} - \lambda_\eps \f{\int x N \ud x}{\int \rho_\alpha * g N_\eps \ud x} | + |\lambda_\eps \f{\int x N \ud x}{\int \rho_\alpha * g N_\eps \ud x} - \lambda_0 \f{\int x N \ud x}{\int \rho_\alpha * g N_\eps \ud x}| \\
\vspace{0.15cm}
& &+\ds|\lambda_0 \f{\int x N \ud x}{\int \rho_\alpha * g N_\eps \ud x} - \lambda_0 \f{\int x N \ud x}{\int g N\ud x} |.
\end{array}$$
The assumptions of Theorem \ref{th:filter} together with Estimates \ref{lm:mollifier:estim4} and \ref{lm:mollifier:estim5} of Lemma \ref{lm:mollifier} give the estimate for $|c_{\eps,\alpha} -c|$ and ends the proof.

\subsection{Quasi-Reversibility Method}

To regularize the exact inverse problem \eqref{statio:inverse}, the so called \emph{quasi-reversibility} method proposed in \cite{P.Z} for the case $\kappa=\delta_{x=\f{y}{2}}$ consisted in adding a term derivative $\alpha\f{\p}{\partial_x} (BN(2x))$ with a small $\alpha >0$ to the right-hand side of Equation \eqref{statio}, viewed as an equation taken in the variable $y=2x.$ The main difference is that we need here to take this term in the variable $x$ and not $2x$ due to the general form of the nonlocal kernel $\kappa.$ We choose to define, for $\alpha >0$ and $k\in \R,$ the following regularised problem
\begin{equation}
\label{(Pr-est)}\left\{\begin{array}{ll}  
{\cal L}_k^\alpha(H_\eps)(x):= \alpha x^{-k}\ds\f{\p}{\partial_x} \bigl(x^{k+1} H_\eps(x)\bigr)
+ H_\eps(x) - 2\int^{\infty}_{0}\kappa(x,y)H_\eps(y) \ud y
= - c_{\alpha,\varepsilon}\ds\f{\p}{\partial_{x}}\bigr(gN_{\varepsilon}(x)\bigl)-\lambda_{\varepsilon}N_{\varepsilon}(x),\\
 \vspace{0.1cm}
H_{\varepsilon}(0)=0; \quad 0<\alpha<1, \quad p\in\mathbb{R}. 
\end{array} \right.\vspace{0.05cm}
\end{equation}
This equation has to be understood in a distribution sense in $\R_+$ undowed with the measure $x^p \ud x.$ We moreover assume that $Supp(N_\eps) \subset \R_+^*.$ Other adaptations would be possible, all consisting in adding a small term derivative of the form $\pm\alpha f_1(x)\f{\p}{\partial_x} (f_2(x) BN (x)),$ with $\alpha >0$ and a boundary condition taken either in $x=0$ if $\alpha>0$ or $x=+\infty$ if $\alpha<0.$ Numerically indeed, $\alpha <0$ proved to give better results (see below Section \ref{sec:num:QR}). The key point is to check that the regularized operator ${\cal L}_k^\alpha$ satisfies Proposition \ref{prop:qr} below.

The choice of $c_{\alpha,\eps}$ is not directly given by integration of the equation, contrarily to the case of \cite{P.Z}. 
Neglecting the regularisation terms involving $\alpha,$ we thus define, as for the exact equation \eqref{statio}:
\begin{equation}\label{def:calpha:qr}
c_{\alpha,\eps}=\ds\f{\lambda_\eps \ds\int x N_\eps (x) \ud x}{\ds\int g(x) N_\eps (x) \ud x}.
\end{equation}

\begin{theor} \label{main_th} 
Let $g,$ $B$ and $\kappa$ satisfy Assumptions \eqref{as:kappa} and \eqref{assum0}--\eqref{assum5}, and moreover $\f{x^\gamma}{g} \in L^2_0$ with $\gamma$ defined in Assumption \eqref{assum4}. Let $(\lambda_0,N)$ the unique eigenpair solution of Problem \eqref{(P)} (as stated in \cite{D.G}).
Let $p>2$ satisfy Assumption \eqref{coer_cond}.
 Let $N_\eps \in L^1\cap L^\infty (\R_+),$ $Supp(N_\eps) \subset \R_+^*,$ and $\lambda_\eps >0$ satisfy $\Vert g(N-N_\eps)\Vert_{L^2(x^p \ud x)} \leq \eps \Vert gN\Vert_{L^2(x^p \ud x)},$ $|\lambda_\eps -\lambda_0|\leq \eps \lambda_0,$ $\Vert N - N_\eps \Vert_{L^1\big((1+x+g(x))\ud x\big)} \leq \eps \Vert N\Vert_{L^1\big((1+x+g(x))\ud x\big)} $ and $\Vert N_\eps - N\Vert_{L^2\big((x^p+1)\ud x\big)} \leq \eps \Vert N\Vert _{L^2\big((x^p+1) \ud x \big)}.$ 
 Let $H_{\eps} \in L^2(x^p \ud x)$ be solution to Problem \eqref{(Pr-est)}\eqref{def:calpha:qr}.

Defining $B_{\eps,\alpha}:=\chi_{N_{\eps}(x) \neq 0} H_{\eps} / N_{\eps}$ we have the following estimates: 
\begin{equation}|c_{\eps,\alpha} - c| \leq C(p,N)   {\eps},\label{est:filter:c}
\end{equation}
\begin{equation}||B_{\eps,\alpha} - B||_{L^2 (x^p N_\eps^2 \ud x)} \leq C(p,N) ( \alpha + \frac{\eps}{\alpha}),\label{est:filter}
\end{equation}
where $C$ is a constant depending on $p,$ $k,$  $\lambda_0,$ $\Vert BN\Vert_{H^1\big((1+x^{p+1})\ud x\big)},$ $\Vert N\Vert_{L^1\big((1+x+g(x))\ud x\big)},$ $\Vert gN\Vert_{L^{1}}$ and $\Vert N\Vert _{H^1\big((x^p+1)\ud x\big)}.$
\end{theor}

\begin{proof}
The estimate for $\vert c_{\eps,\alpha} -c\vert$ is obtained in a similar manner as for the filtering method.

For the estimate for $B,$ we first write
$$\Vert B_{\eps,\alpha} N_\eps - BN_\eps\Vert_{L^2(x^p \ud x)} \leq \Vert B_{\eps,\alpha} N_\eps - BN\Vert_{L^2(x^p\ud x)} + \Vert BN - BN_\eps\Vert_{L^2(x^p\ud x)}.$$
The second term of the right-hand side is simply bounded by $\Vert B\Vert_{L^2(x^p\ud x)}\Vert N-N_\eps\Vert_{L^2(x^p\ud x)}\leq \eps   \Vert B\Vert_{L^2(x^p\ud x)}.$ 

For the first term of the right-hand side, as for the filtering method, we decompose $H_\eps -BN,$ and for this we need to establish some regularity properties of the operator ${\cal L}_k^\alpha$ defined in Equation \eqref{(Pr-est)} and designed to approximate ${\cal L}_\kappa.$ This is given by the following proposition, which is for the quasi-reversibility method the equivalent of Lemma \ref{lm:mollifier} for the filtering method.
\begin{prop}\label{prop:qr}
Let $p>2,$ $F=f_1 + \f{\p}{\partial_x} f_2$ with $f_1 \in L^1\big((1+x)\ud x\big)\cap L^2(x^p\ud x)$ and $f_2\in H^1((1+x^p) \ud x)\cap W^{1,1}(x\ud x).$ Let $\kappa,g,p$ satisfy the assumptions of Theorem \ref{main_th}. There exists $u \in L^1(x\ud x)$ solution of the following problem, where $k\neq -2$ and $0<\alpha \leq 1:$
\begin{equation}\label{eq:qr:probgen}
{\cal L}_k^\alpha (u)(x):= \alpha x^{-k}\f{\p}{\partial_x} \bigl(x^{k+1} u\bigr) + {\cal L}_\kappa (u)
=F.
\end{equation}
Moreover, we have the following estimates for a constant $C>0$ only depending on $g,$ $\kappa,$ $k$ and $p:$
\begin{enumerate}
\item $\Vert u \Vert_{L^2(x^p \ud x)} \leq \ds\f{1}{1-2\sqrt{C_p}} \Vert F\Vert_{L^2(x^p \ud x)},$ \label{qr:estim1}\\
\item $\Vert u \Vert_{L^2(x^p \ud x)} \leq \ds\f{C}{\alpha} \Vert f_1+f_2(1+\f{1}{x^2}) \Vert_{L^2(x^p \ud x)}.$\label{qr:estim2}
\end{enumerate}
\end{prop}
\begin{proof}
Let us first establish the existence of a solution in $L^{1}(x \ud x).$ We rewrite \eqref{eq:qr:probgen} as follows
\begin{equation}\label{eq_proof}
 \left\{\begin{array}{ll}
 \alpha x\ds\f{\p}{\partial_{x}}(u(x))+(\alpha (k+1)+1)u(x)=2\ds\int_{x}^{\infty}\kappa(x,y)u(y)\ud y +F(x). \\
 \vspace{0.05cm}
 u(0)=0, \quad p>2.
 \end{array} \right.\vspace{0.05cm}
\end{equation}
We consider $v\in L^{1}(\mathbb{R_{+}},x \ud x)$ and define $u=T(v)$ the explicit solution of  
 \begin{equation*}
 \left\{\begin{array}{ll}
 \alpha x\ds\f{\p}{\partial_{x}}(u(x))+(\alpha (k+1)+1)u(x)=2\ds\int_{0}^{\infty}v(y)\kappa(x,y)\ud y -F(x), \\
 \vspace{0.05cm}
 u(0)=0,\quad p>2.
 \end{array} \right.\vspace{0.05cm}
\end{equation*}
Let $v_{1}$ and $v_{2}$ two functions of $L^{1}(\mathbb{R_{+}},x \ud x)$ associated to $u_{1}$ and $u_{2}$ then by doing the difference between the two equations satisfied in the one hand by $u_{1}, v_{1}$ and in the other hand by $u_{2}, v_{2}$ we have
\begin{equation}\label{eq_u}
\alpha x\f{\p}{\partial_{x}}(\delta u(x))+(\alpha (k+1)+1)\delta u(x)=2\ds\int_{0}^{\infty}\delta v(y)\kappa(x,y)\ud y, \quad \text{with} \; \delta u=u_{1}-u_{2} \;\text{and}\;  \delta v=v_{1}-v_{2}, 
\end{equation}
what implies the inequality (see \cite{P}, prop.6.3 for instance)
 $$\alpha x\f{\p}{\partial_{x}}\vert \delta u(x)\vert+(\alpha (k+1)+1)\vert \delta u(x)\vert\leq 2\ds\int_{0}^{\infty}\vert \delta v(y)\vert \kappa(x,y)\ud y.$$
 Multiplying by $x$ and integrating on $[0,\infty[$ we deduce the estimate
 $$\int_{0}^{\infty}x\vert \delta u(x)\vert \ud x \leq \frac{1}{\alpha(k-1)+1}\int_{0}^{\infty}y\vert \delta v(y)\vert \ud y.$$
 This proves that $T$ is a Lipschitz function and we deduce the existence of a solution $u\in L^{1}(\mathbb{R}_{+},x\ud x)$ by the Schauder fixed point theorem.\\
 
 For the first estimate, we multiply Equation \eqref{eq:qr:probgen} by $x^p u$ and integrate from $0$ to $x.$ Using that $u x^{p-k} \ds\f{\p}{\partial_x} (x^{k+1} u) = (k+1)x^p u^2 + x^{p+1} \ds\f{\p}{\partial_x} (\f{u^2}{2}),$ it gives
 $$\int_0^x \alpha (k+1) x^{p} u^2(x) \ud x  + \f{\alpha}{2} x^{p+1} u^2(x) + \int_0^x {\cal L}_\kappa (u)(x) u (x) x^p \ud x = \int_0^x F(x) u(x) x^p \ud x.$$
 From this, we deduce 
  \begin{equation}\label{estim:qr}\int {\cal L}_\kappa (u) (x) u(x) x^p \ud x \leq \int  F(x) u(x) x^p \ud x.\end{equation}
 Applying the coercitivity on $L^2(x^p \ud x)$ of the bilinear form ${\cal A} (u,v)=\ds\int {\cal L}_\kappa (u) v x^p \ud x$ we get immediately the first estimate.
 
 For the second one, we integrate by part, on the right-hand side of Equation \eqref{eq:qr:probgen}, the term with $\f{\p}{\partial_x} f_2,$ and use the equation to express $\f{\p}{\partial_x} ( u)$ with the other terms of the equation:
 $$\begin{array}{lll}
 \ds\int (\f{\p}{\partial_x} f_2)ux^p \ud x &=&-\ds\int f_2 \f{\p}{\partial_x} (x^pu) \ud x = -\ds\int \f{p }{x} f_2  u x^p \ud x -  \int f_2 x^p \f{\p}{\partial_x} (u) \ud x\\ \\ &=&-\ds\int \f{p }{x} f_2  u x^p \ud x +  \ds\int \f{k+1}{x} f_2  u x^p \ud x  
 + \ds\int \f{1}{\alpha x}  f_2 \big({\cal L}_\kappa (u)(x) - f_1(x) - \f{\p}{\partial_x} f_2 \big)x^p\ud x \\ \\
&=&\ds\int \f{k+1-p }{x} f_2  u x^p \ud x   
 + \ds\f{1}{\alpha}\biggl(\int \f{1}{x} f_2 \big({\cal L}_\kappa (u)(x) - f_1(x)\big) x^p \ud x + \ds\f{p-1}{2} \int f_2 ^2 x^{p-2} \ud x \biggr)  \\ \\
 &\leq & \ds\f{C}{\alpha} \Vert f_2\Vert_{L^2((x^{p-1}+x^p)\ud x)} \big(\Vert u\Vert_{L^2(x^p\ud x)} + \Vert f_1\Vert_{L^2(x^p\ud x)}\big) + \f{p}{\alpha}\Vert f_2\Vert_{L^2(x^{p-2}\ud x)}^2 \end{array}.$$
 Together with the first estimate, it provides the desired inequality.
\end{proof}

We are now ready for the proof of Theorem \ref{main_th}. We see that $H_\eps$ can be viewed as solution of Equation \eqref{eq:qr:probgen} with $F=- c_{\alpha,\varepsilon}\ds\f{\p}{\partial_{x}}\bigr(gN_{\varepsilon}(x)\bigl)-\lambda_{\varepsilon}N_{\varepsilon}(x),$ whereas $H=BN$ would be solution of \eqref{eq:qr:probgen} if $\alpha=0$ and 
$F=- c\f{\p}{\partial_{x}}\bigr(gN(x)\bigl)-\lambda_0 N(x).$
To isolate in the error term the contribution due to the $\alpha-$regularization from the one due to the measurement error $\eps,$ we thus define an intermediate function $H_\alpha$ as the solution of Equation \eqref{eq:qr:probgen} with $F=- c\ds\f{\p}{\partial_{x}}\bigr(gN(x)\bigl)-\lambda_0 N(x).$ We then write:
$$\Vert B_{\eps,\alpha} N_\eps - BN\Vert_{L^2(x^p\ud x)} = \Vert H_\eps - H\Vert_{L^2(x^p\ud x)}  \leq \Vert H_\eps - H_\alpha\Vert_{L^2(x^p\ud x)}  + \Vert H_\alpha - H\Vert_{L^2(x^p\ud x)}. 
$$
The function $H_\eps - H_\alpha$ is solution of Equation \eqref{eq:qr:probgen} with 
$$F_{\eps,\alpha}=\f{\p}{\partial_{x}}\bigr(- c_{\alpha,\varepsilon} gN_{\varepsilon}(x)+ c gN(x)\bigl)
-\lambda_{\varepsilon}N_{\varepsilon}(x) +\lambda_0 N(x),$$
and we can use Estimate \ref{qr:estim2} of Proposition \ref{prop:qr} to obtain an error term in the order of $\ds\f{\eps}{\alpha}.$ The difference $H_\alpha - H$ is solution of Equation \eqref{eq:qr:probgen} with $F=-\alpha x^{-k}\ds\f{\p}{\partial_x} (x^{k+1} BN),$ and we can use Estimate \ref{qr:estim1} of Proposition \ref{prop:qr} to obtain an error term bounded by $C \alpha \Vert BN\Vert_{H^1\big((1+x^{p+1})\ud x\big)}.$ It ends the proof of Theorem \ref{main_th}.

\end{proof}

\section{Numerical approach of the inverse problem}
\subsection{The direct problem}\label{direct}
 Assuming that the division rate $B,$ the growth rate $g$ and $c>0$ are known, we solve the time-dependent problem \eqref{P_t} and look for a steady dynamics.\\
 We choose to split the time evolution of the problem into its conservative advection part and into its gain and lost part by division as follows 
 \begin{equation*}
 \left\{\begin{array}{ll}  
 \vspace{0.15cm}
 \ds\f{\p}{\partial_{t}}n(t,x)+c\f{\p}{\partial_{x}}(g(x)n(t,x))=0 \\
 \ds\f{\p}{\partial_{t}}n(t,x)+B(x)n(t,x)=2\ds\int_{x}^{\infty}B(y)\kappa(x,y)n(t,y)\ud y.\\
 \end{array} \right.\vspace{0.05cm}
\end{equation*}
 We use an upwind finite volume method with computation length domain $L$ and grid number points $ka$: $x_{i}=i\Delta x, \; \,  0\leq i\leq ka$ with $\Delta x=L/ka$ \\
 $$n_{i}^{k}=\frac{1}{\Delta x}\int_{x_{i-\frac{1}{2}}}^{x_{i+\frac{1}{2}}} n(k\Delta t,y)\ud y, \quad \frac{1}{\Delta t}\int_{0}^{\Delta t} n(k\Delta t +s,x_{i})\ud s \thickapprox n_{i}^{k+1}. $$
For the time discretization one can choose, thanks to the CFL ({\it{Courant-Friedrichs-Lewy}}) stability condition, the time step $\Delta t <\ds\frac{1}{\ds\max_{i\in{1, ...,ka}}(B_{i}+\frac{c}{\Delta x}g_{i})}$ with the notation $g_{i}=g(i\Delta x)$ and $B_{i}=B(i\Delta x)$. \\
The numerical scheme is given for $i=1, ...,ka$ by $n_{0}^{k}$ and 
\begin{itemize}
\item[\textbullet] First for the conservative equation
$$n^{k+1/2}_{i}=n^{k}_{i}-c\frac{\Delta t}{\Delta x}\bigl((gn)^{k}_{i+1/2}-(gn)^{k}_{i-1/2}\bigr),$$
the interface fluxes $(gn)^{k}_{i\pm1/2}$ are defined by upwind method. \\ 

\item[\textbullet] Second for the gain and loss part by cellular division we compute 
$$n^{k+1}_{i}=\bigl(1-\Delta t\, B_{i}\bigr)n^{k+1/2}_{i}+2\Delta t \,\mathcal{F}^{k}_{i} $$   
where $\mathcal{F}^{k}_{i} \thickapprox \ds\int_{x_{i}}^{x_{ka}} B(y)n^{k+1/2}(y)\kappa(x_{i},y)\ud y$.\\
\item[\textbullet] At last we renormalize the discrete solution by
$$\tilde{n}^{k+1}=\frac{n^{k+1}}{\ds \sum_{j=1}^{ka}n^{k+1}_{j}\Delta x}$$
what allows to have $\tilde{n}^{k+1}\xrightarrow[k\to\infty]{}N$, \quad  $\ds \sum_{i=1}^{ka}N_{i}\Delta x=1$, \quad $N_{i}>0,$ where $N$ is the dominant eigenvector for the discrete problem associated to the following steady equation
$$ c\f{\p}{\partial_{x}}(g(x)N(x))+(B(x)+\lambda_{0})N(x)=2\ds\int^{+\infty}_{0}B(y)\kappa(x,y)N(y)\ud y$$ with $\lambda_{0}$ the dominant eigenvalue associated to $N$.
\end{itemize}  
\subsection{The inverse problem without regularization}
As illustrated in \cite{D.P.Z}, solving numerically  Equation \eqref{statio:inverse} without regularization is unstable. Indeed, this recovering naive method gives bad reconstructions of $H=BN$ as soon as the observed $N_\eps$ is irregular (see above the estimates and see also \cite{D.P.Z}). Here, what we moreover observe is that, at the neighborhood of $x=0$, the solution explodes. As an example of this we consider a length domain $L=25$ and the total number grid points $ka=300$. By an upwind method we compute numerically the Equation \eqref{statio:inverse} and compare the result with that obtained by solving the direct problem (D.P).

\begin{figure}[!ht]
\resizebox{15cm}{8.5cm} {\input{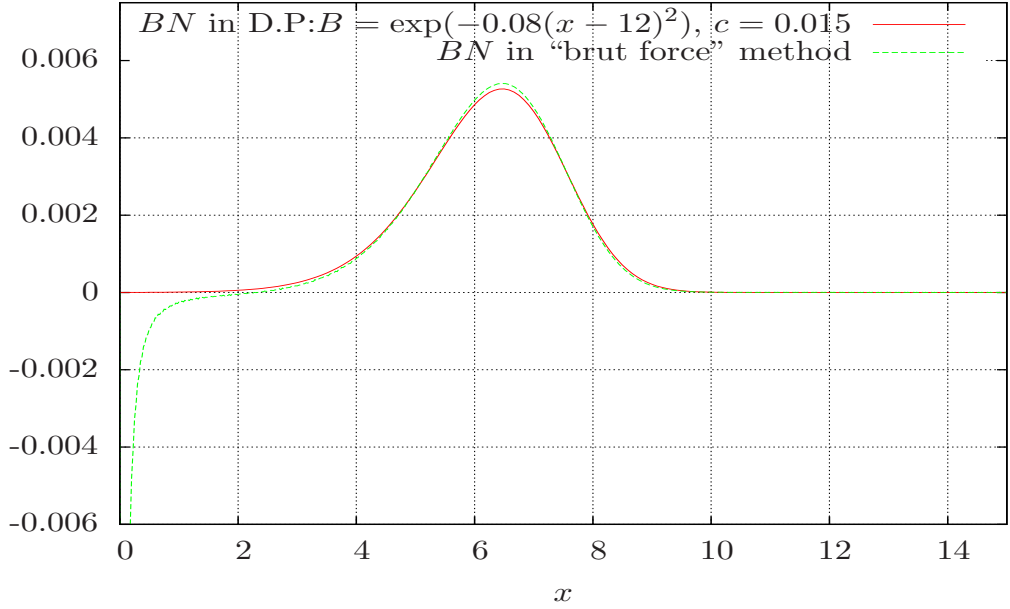}}
\caption{Numerical reconstruction of $BN$ by ``brute force method" with the choice $g(x)=x$ and $\kappa(x,y)=\ds\frac{1}{y}\mathbb{I}_{x<y}$.}
\end{figure}     

\subsection{The inverse problem : Quasi-Reversibility discretization} \label{sec:num:QR}
In this section we numerically investigate the regularization of the inverse problem \eqref{statio:inverse} by the Quasi-reversibility method. It is based on Equation \eqref{(Pr-est)}, that we rewrite, dropping the index $\varepsilon$, as follows 
\begin{equation*}
\left\{\begin{array}{ll}  
\vspace{0.15cm}
 -\alpha x^{-k}\ds\f{\p}{\partial_x} (x^{k+1}B_{\alpha}(x)N(x))+c\ds\f{\p}{\partial_{x}}(g(x)N(x) )+(B_{\alpha}(x)+\lambda_{0})N(x)  =2\ds\int^{\infty}_{x}\kappa(x,y)B_{\alpha}(y)N(y) \ud y,\\
(B_{\alpha}N)(0)=0, \quad (B_{\alpha}N)(\infty)=0; \quad 0<\alpha; \; \;  k\in \R_{+}.
\end{array} \right.\vspace{0.05cm}
\end{equation*}
Assuming that $N$ and $\lambda$ are measured, we first define $c$ by \eqref{def:calpha:qr} and then look for an estimate of the division rate $B_{\alpha}$. For this, 
  we put the notation $$H_{\alpha}=B_{\alpha}N\;  \text{and} \; \; L=-c\ds\f{\p}{\partial_{x}}(gN)-\lambda_{0}N.$$
By a standard upwind method we obtain, when dropping the index $\alpha$, the following discretization
 \begin{equation*}
\left\{\begin{array}{ll}  
\ds -\alpha x_{i}^{-k}\biggl(\frac{x_{i+1}^{k+1}H_{i+1}-x_{i}^{k+1}H_{i}}{\Delta x}\biggr)+H_{i}-2\ds\sum_{j=i}^{ka}H_{j}\kappa_{i,j}\Delta x =L_{i} \\
 \vspace{0.1cm}
\text{with} \; \;  L_{i}=-\lambda_{0}N_{i} -c\biggl(\ds\frac{g_{i+1}N_{i+1}-g_{i}N_{i}}{\Delta x}\biggr), \quad \forall \; i=1,...,ka \\
H_{0}=0 \;  \text{and} \;  \; H_{l}=0, \quad \forall \; l>ka. \\ 
 \end{array} \right.\vspace{0.05cm}
\end{equation*}
By developing this discrete equation we obtain
$$\biggl(-\alpha\frac{(i+1)^{k+1}}{i^k} -2\kappa_{i,i+1}\Delta x\biggr)H_{i+1}+\biggl(1+\alpha i -2\kappa_{i,i}\Delta x\biggr)H_{i}-2\ds\sum_{j=i+2}^{ka}H_{j}\kappa_{i,j}\Delta x =L_{i}, \quad \forall \; i=1,...,ka.$$
We rewrite it under matrix shape  
$A\times H=L$ with $A$ a matrix of coefficients of size $ka\times ka$; $H$ is the unknown vector of size $ka$ and $L$ is a known vector of size $ka$.\\
The matrix $A$ being a upper triangular one, we can
 solve directly the linear system thanks to the following iterations
\begin{equation*}
\left\{\begin{array}{ll}  
H_{ka}=\ds\frac{L_{ka}}{A_{ka,ka}}; \\
H_{i}=\ds\frac{1}{A_{i,i}}\biggl(L_{i}-\ds\sum_{j=i+1}^{ka}A_{i,j}H_{j}\biggr) \,; \quad \forall \; i=ka-1,...,1 
 \end{array} \right.\vspace{0.05cm}
\end{equation*}
The matrix $A$ satisfying $A_{ij} >0$ for $j\geq i+1,$ we can choose $\Delta x$ small enough so that $A_{ii}=1+\alpha i - 2 \kappa_{i,i} \Delta x >0$ for all $i.$ This guarantees that no oscillations appear. 

\subsection{The inverse problem: Filtering discretization} 
This section is devoted to the numerical discretization of the inverse problem \eqref{statio:inverse} by the Filtering method based on Equations \eqref{Pestf}-\eqref{def:mollifier}. 
The aim is to numerically solve the Equation \eqref{Pestf} that we rewrite when dropping the index $\varepsilon$ as follows
\begin{equation*}
\left\{\begin{array}{ll}  
\vspace{0.15cm} 
 x^{p}c_{\alpha}\ds\f{\p}{\partial_{x}}\bigl(g(x)N_{\alpha}(x)\bigr)+ x^{p}(B_{\alpha}(x)+\lambda_{\alpha})N_{\alpha}(x)=2 x^{p}\ds\int^{+\infty}_{0}\kappa(x,y)B_{\alpha}(y)N_{\alpha}(y)\ud y, \\
\bigl(B_{\alpha}N_{\alpha}\bigr)(0)=0 , \quad \alpha>0,
\end{array} \right.\vspace{0.05cm}
\end{equation*}
with $N_{\alpha}=N$\,\textasteriskcentered \,$\rho_{\alpha}$ and $\rho_{\alpha}$ a sequence of mollifiers.\\
As previously, we want to  estimate $B_{\alpha}$ from a measured density $N$ and Malthus parameter $\lambda.$ We first define $c$ by \eqref{def:cepsalpha}.
We then rewrite the regularised equation as follows
$$ x^{p}B_{\alpha}(x)N_{\alpha}(x)-2\ds\int_{x}^{\infty}x^{p}\kappa(x,y)B(y)N_{\alpha}(y)\ud y = -x^{p}c_{\alpha}\ds\f{\p}{\partial_{x}}\bigl(g(x)N_{\alpha}(x)\bigr)-\lambda_{\alpha} x^{p}N_{\alpha}(x)$$
For the convolution terms arising in the previous equation we use the combination of the Fast Fourier Transform and its inverse which we respectively note by $F$ and $F^{*}$ then we define the mollifiers $\rho_{\alpha}$ by its Fourier transform: $\hat{\rho}_{\alpha}(\xi)=\ds\frac{1}{\sqrt{1+\alpha^{2}\xi^{2}}}$.\\
This leads to the following approximations
$$N_{\alpha}\thickapprox F^{*}\bigl(\hat{\rho}_{\alpha}(\xi)F(N)(\xi)\bigr); \quad \f{\p}{\partial_{x}}\bigl(g N_{\alpha}\bigr)\thickapprox dGN_{\alpha}= F^{*}\bigl(i\xi\hat{\rho}_{\alpha}(\xi)F(gN)(\xi)\bigr). $$ 
 For the discretization we put the notation $$H_{\alpha}=B_{\alpha}N_{\alpha}\,\text{and} \quad 
 L_\alpha= -c_{\alpha}dGN_{\alpha}-\lambda_{\alpha}N_{\alpha}$$ 
 then in each grid point $x_{i}=i\Delta x$ we obtain when dropping the index $\alpha$:
\begin{equation*}
\left\{\begin{array}{ll}  
H_{0}=0 \\
x_{i}^{p}(1-2\kappa_{i,i}\Delta x)H_{i}-2\ds\sum_{j=i+1}^{ka}x_{i}^{p}H_{j}\kappa_{i,j}\Delta x =x_{i}^{p}L_{i}\, ; \quad \forall \; i=1,...,ka.
\end{array} \right.\vspace{0.05cm}
\end{equation*} 
 We rewrite this previous discrete equation under matrix shape $A\times H=L$ with $A$ the matrix of coeficients which is an upper triangular one and of size $ka\times ka$.\\
 The shape of the matrix $A$ allows to use adequately the LU iterative numerical method, and then we deduce the following iteration
 \begin{equation*}
\left\{\begin{array}{ll}  
H_{ka}=\ds\frac{L_{ka}}{A_{ka,ka}} \\
H_{i}=\ds\frac{1}{A_{i,i}}\biggl(L_{i}-\ds\sum_{j=i+1}^{ka}A_{i,j}H_{j}\biggr) \,; \quad \forall \; i=ka-1,...,1 
\end{array} \right.\vspace{0.05cm}
\end{equation*} 

\section{Numerical Tests}
For the numerical tests we use as input data the noisy one $N_{\varepsilon}$ which correspond for $\varepsilon=0$ to the eigenfunction $N$ obtained by solving numerically the long time behavior of the direct problem in section \ref{direct}. 
The direct problem is solved in the length domain $L=25$ for $ka=300$ number grid points with two differents initial data: a step initial data and a maxwellian one, as follows
\begin{equation}\label{init}
\text{Step function}:\left\{\begin{array}{ll}  
n^{0}(x)=0.2 \quad 5\leq x\leq 10 ,\\
n^{0}(x)=0 \quad \text{other where} \\
\end{array} \right.\vspace{0.05cm}
\; \hfill
\text{Maxwellian}:\left\{\begin{array}{ll}  
n^{0}(x)= \ds\frac{1}{\sqrt{0.4\pi}}\exp\biggl(-\frac{(x-10)^{2}}{0.4}\biggl)   ,\\
\forall \; 0\leq x\leq L
\end{array} \right.\vspace{0.05cm}
\end{equation} 
and the steady solution is taken when $\Vert n(t,x)-N(x)\Vert_{L^{1}}<=10^{-10}$.\\

\begin{figure}[!ht]
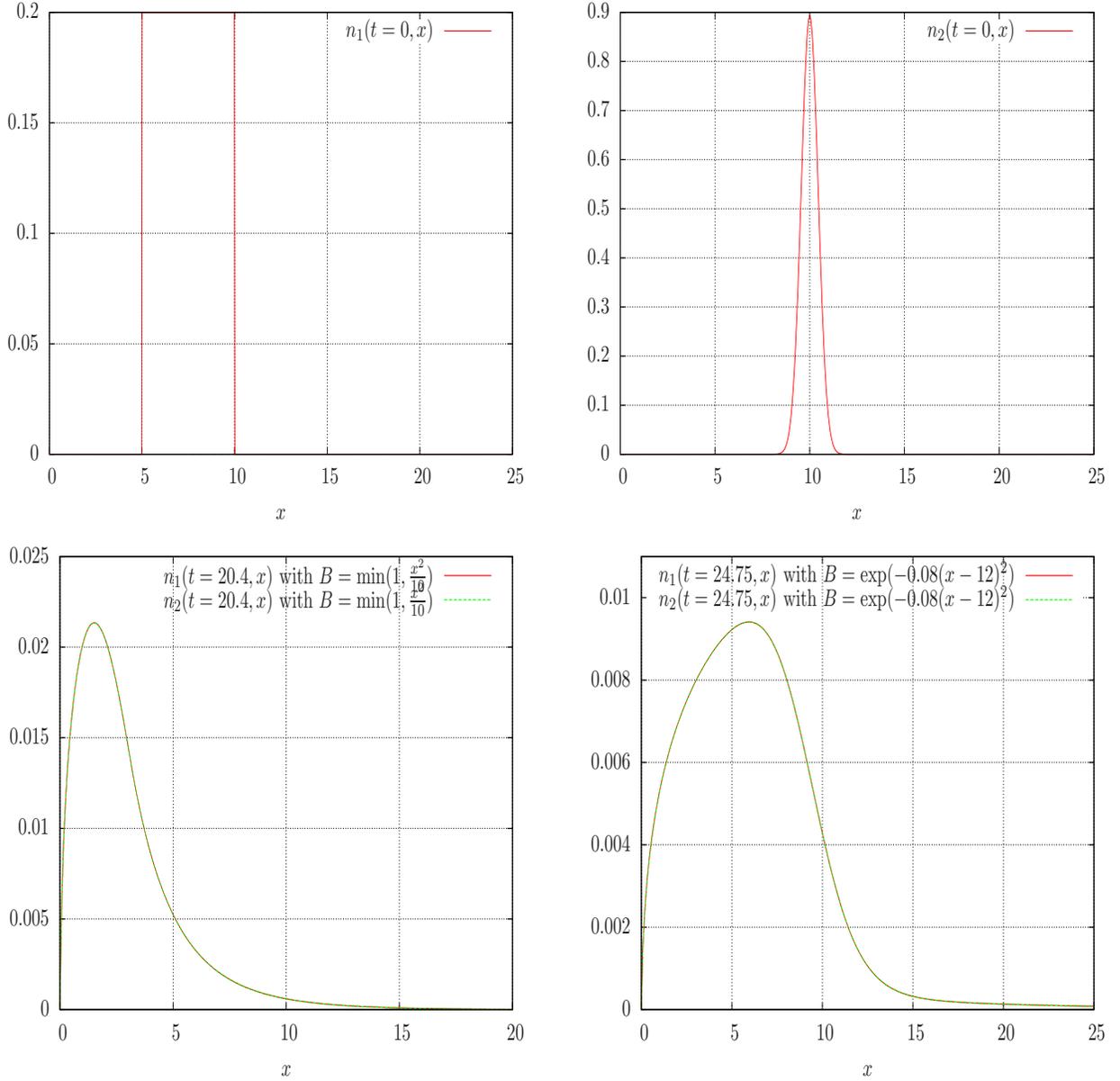

\resizebox{8.5cm}{8cm} {\input{step_initial.tex}} \hfill
\resizebox{8.5cm}{8cm} {\input{maxw_initial.tex}} \vfill
\resizebox{8.5cm}{8.2cm} {\input{profil_asym1.tex}} \hfill
\resizebox{8.5cm}{8.2cm} {\input{profil_asym2.tex}}
\caption{Direct problem $g=x^{1/2}$: Top left: Step initial function. Top right :Maxwellian initial function. Down left: Steady solutions of cellular density with $c=1$. Down right: Steady solutions of cellular density with $c=0.5$ .}
\label{fig1}
\end{figure}

In order to show the unique asymptotic profile of the direct problem we plot in pictures Fig.\ref{fig1} the steady cellular density $N$ related to the two previous initial data with different values of $c$, $B$ and with the choice $g(x)= x^{1/2}$ and $\kappa(x,y)=\ds\frac{1}{y}\,\mathbb{I}_{\{x<y\}}$. 
\newpage

\subsection{Numerical reconstruction of $BN$ in the noiseless case $\varepsilon=0$} 
For the case where the input data are exactly known i.e. for $\varepsilon=0$, we recover thanks to the Quasi-reversibility and Filtering methods the division rate $B$ by computing numerically the value of $BN$ with $N$ obtained by solving the direct problem with high precision and for various choices of the division rate $B$ as shows in figure Fig.\ref{allure_B} below. \\
\begin{figure}[!ht]
\resizebox{16cm}{9.8cm} {\input{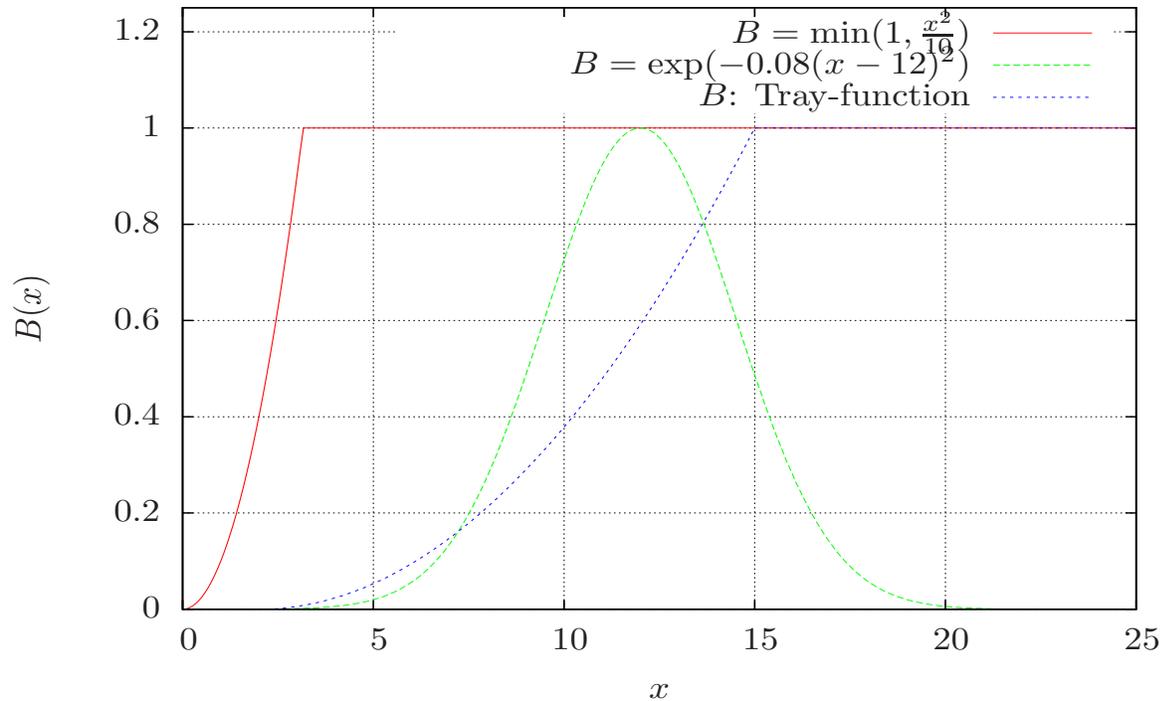}}
 \caption{Various choices of $B$ to solve the direct problem .}
\label{allure_B}
\end{figure} 
\\
In figure Fig.\ref{allure_B}, the Tray-function is defined as follow
\begin{equation*}
B(x)=\ds\left \{  \begin{array}{ll} 
0, \quad \text{for} \; \; x<2 \\
\\
\ds\frac{(x-2)^2} {13^2}, \quad \text{for} \; \;  x\in [2,15] \\
\\
1,  \quad \text{for} \; \;  x>15 .\\
\end{array} \right.
\end{equation*}
Then with the notation D.P for the direct problem we obtain: \\
\begin{itemize}
\item[\textbullet] {\bf{For the choice $\kappa(x,y)=\ds\frac{1}{y}\,\mathbb{I}_{\{x<y\}}$}} 
\end{itemize}
  \begin{figure}[!ht]
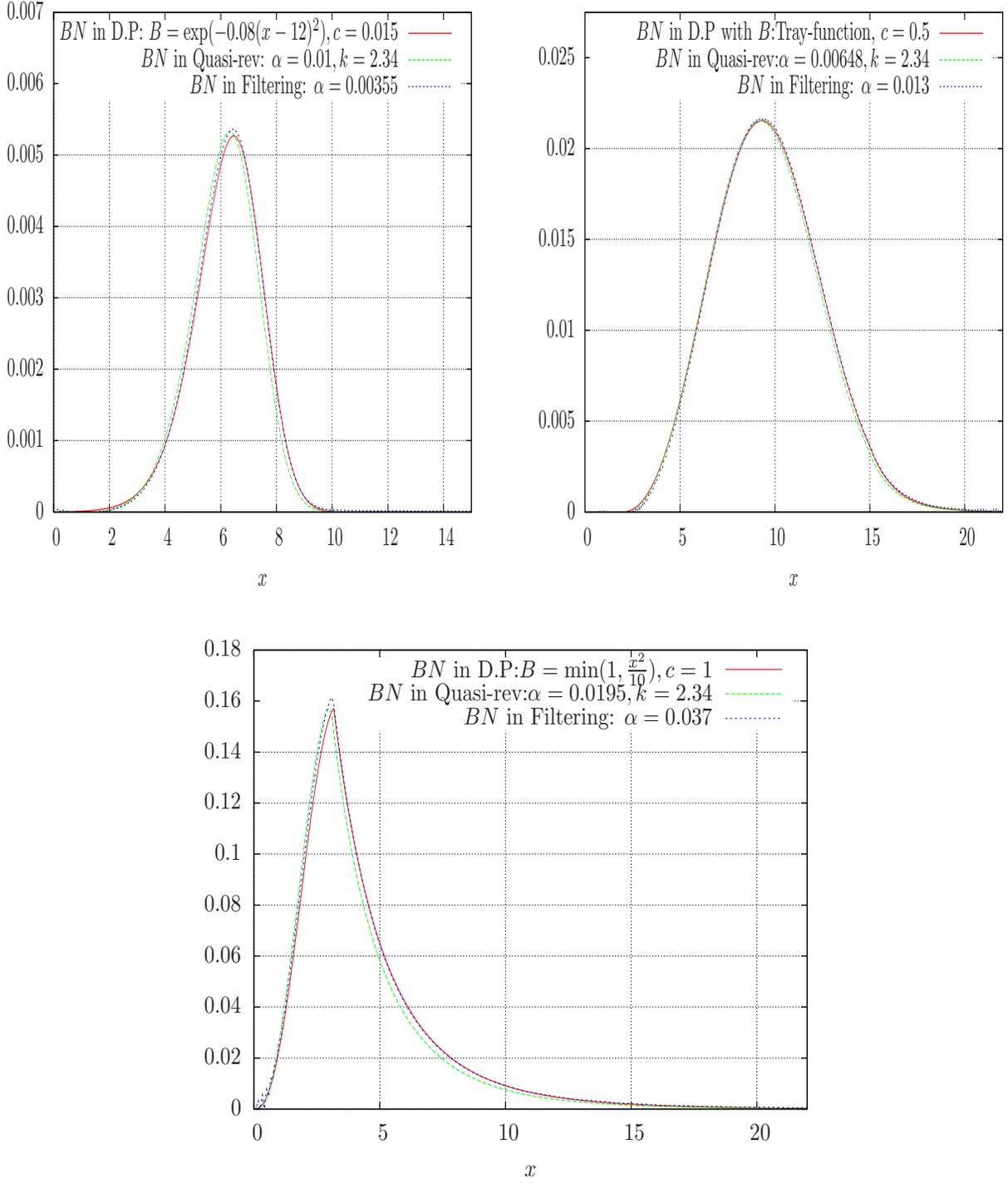

\resizebox{8.5cm}{9.8cm} {\input{comp1.tex}}\hfill
 \resizebox{8.5cm}{9.8cm} {\input{comp2.tex}}\vfill
 \begin{center}
 \resizebox{11cm}{9cm} {\input{comp3.tex}}
 \end{center}
 \caption{Numerical reconstruction of $BN$ for each regularization method in the case $\varepsilon=0$. \qquad \qquad Top left :$g(x)=x$. Top right :$g(x)= x^{1/3}$. Down :$g(x)=x^{1/2}$ . }
\label{fig3}
\end{figure} 
 \newpage
 We measure the relative error in $L^2$ norm by
 \begin{equation}\label{eq_err}
 error=\ds\frac{\bigl\Vert BN-(BN)_{\varepsilon,\alpha}\bigr\Vert_{L^2}}{\bigl\Vert BN\bigr\Vert_{L^2}}, 
 \end{equation}
 where $BN$ is the exact numerical solution of the direct problem and $(BN)_{\varepsilon,\alpha}$ represents the numerical reconstruction either by the Quasi-reversibility method or by  the Filtering one. So we obtain for instance for the given parameters $g(x)=x, \; c=0.015$ and $k=2.34$ the following reconstruction error of the division rate as a function of $\alpha$. 
 \begin{figure}[!ht]
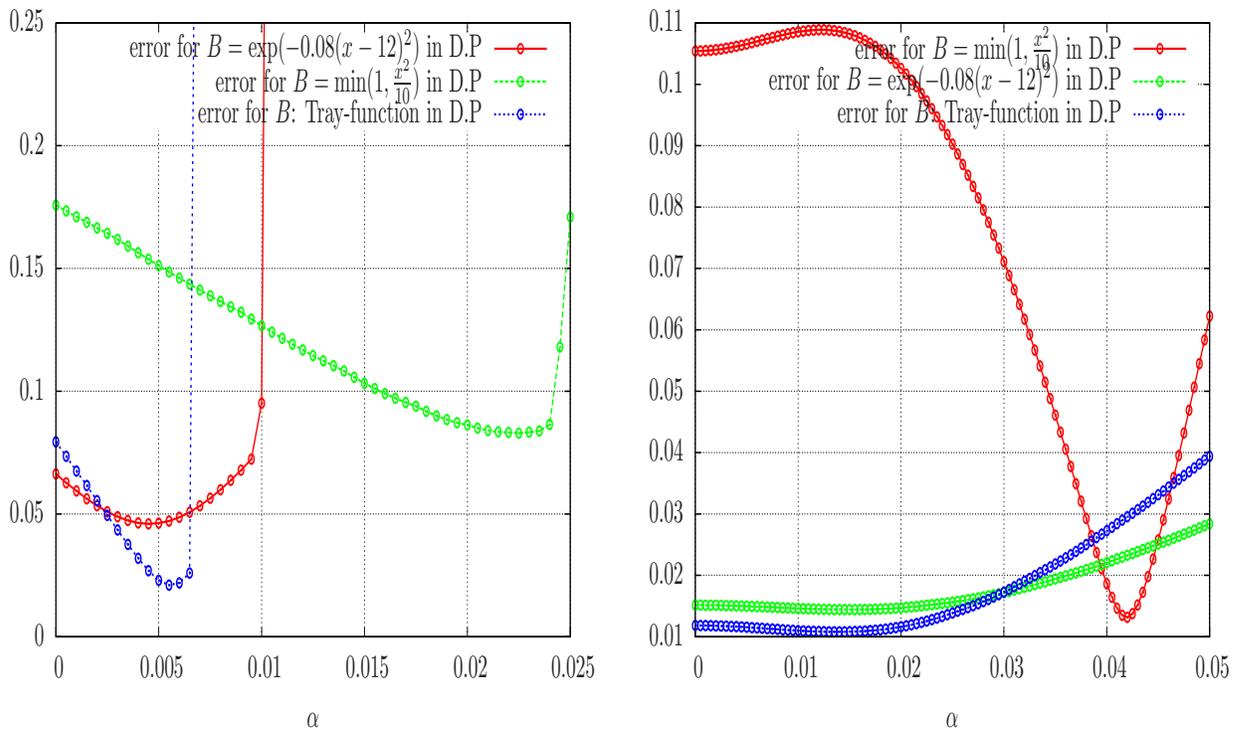

\resizebox{8.5cm}{10cm} {\input{leon_marie_1.tex}}\hfill
\resizebox{8.5cm}{10cm} {\input{leon_marie_2.tex}}
 \caption{Numerical errors for $\varepsilon=0$ with different choices of $B$ in the direct problem. 
 Left : errors by Quasi-reversibility method. Right : errors by Filtering method.}
\label{fig4}
\end{figure} 
\\
  \begin{itemize}
\item[\textbullet] {\bf{For the choice $\kappa(x,y)=\ds\frac{1}{y}\kappa_{0}(\frac{x}{y})$ with $\kappa_{0}\sim\mathcal{N}(\frac{1}{2},\frac{1}{4})$}} 
 \begin{figure}[!ht]
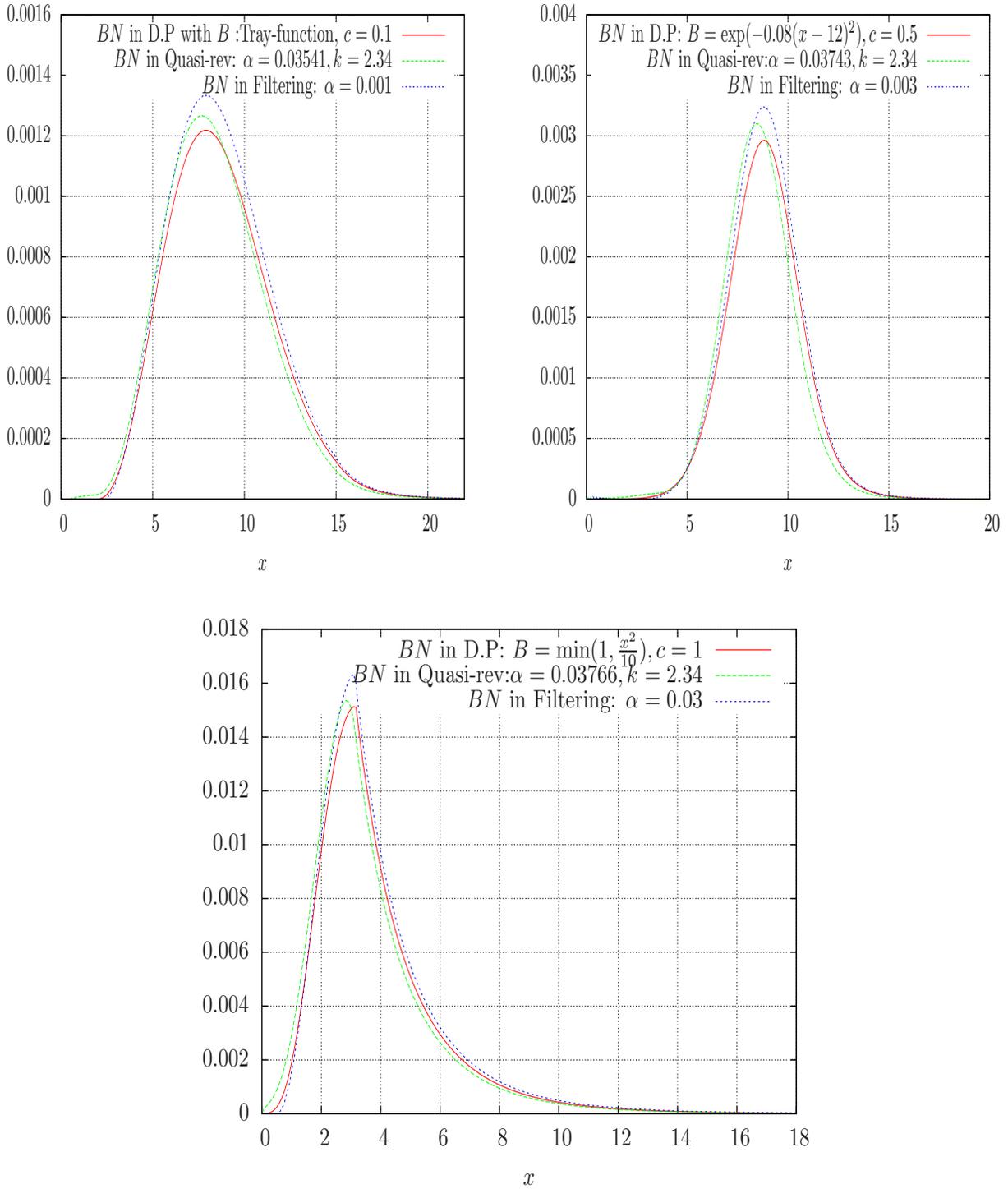

\resizebox{8.5cm}{9.6cm} {\input{comp4.tex}}\hfill
\resizebox{8.5cm}{9.6cm} {\input{comp5.tex}}
 \begin{center}
 \resizebox{11cm}{9.6cm} {\input{comp6.tex}}
 \end{center}
\caption{Numerical reconstruction of $BN$ for each regularization method in the case $\varepsilon=0$. \qquad \qquad Top left:$g(x)= x$. Top right: $g(x)=x^{1/3}$. Down : $g=x^{1/2}$ .}
\label{fig5}
\end{figure} 
\newpage
We measure the reconstruction error thanks to the relation \eqref{eq_err} for the given parameters $g(x)= x^{1/3}, \; c=0.5,\; k=2.34$ and we obtain the following representations
 \begin{figure}[!ht]
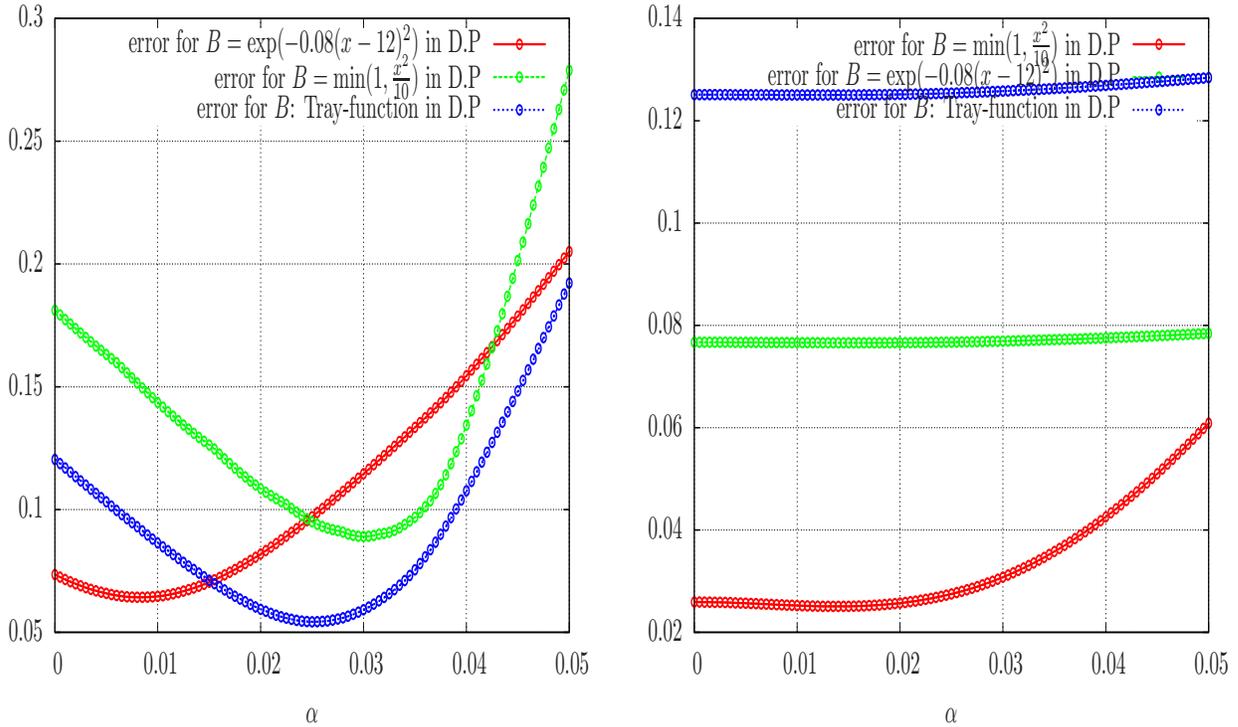

\resizebox{8.5cm}{10cm} {\input{leon_marie4.tex}}\hfill
\resizebox{8.5cm}{10cm} {\input{leon_marie_3.tex}}
 \caption{Numerical errors for $\varepsilon=0$ with different choices of $B$ and $c$ in the direct problem. \qquad \qquad 
 Left : errors by Quasi-reversibility method ($k=2.34$) . Right : errors by Filtering method.}
\label{fig6}
\end{figure} 
\\
\end{itemize}
\subsection{Numerical reconstruction of $BN$ in the noisy case $\varepsilon\neq0$} 
For this case, we consider as input data the values of the solution $N$ of the direct problem in which we add a multiplicative random noise uniformely distributed in $[\frac{-\varepsilon}{2},\frac{\varepsilon}{2}]$ (see \cite{DHRR} for a more precise statistical setting of noisy informations). The nonnegativity of the data is insured by the choice $$N_{\varepsilon}=\max(N(1+l\varepsilon),0), \quad l\in [-\frac{1}{2},\frac{1}{2}], \; \varepsilon\in[0,1]. $$ 
Then with these noisy data we numerically obtain\\
\begin{itemize}
\item[\textbullet] {\bf{For the case $\kappa(x,y)=\ds\frac{1}{y}\mathbb{I}_{\{x<y\}}$}}
 \begin{figure}[!ht]
\resizebox{8.5cm}{9.75cm} {\input{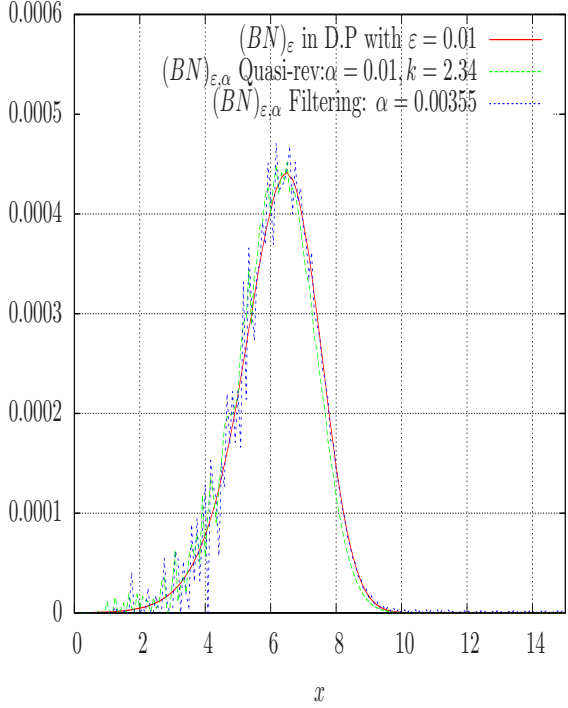}} \hfill
\resizebox{8.5cm}{9.75cm} {\input{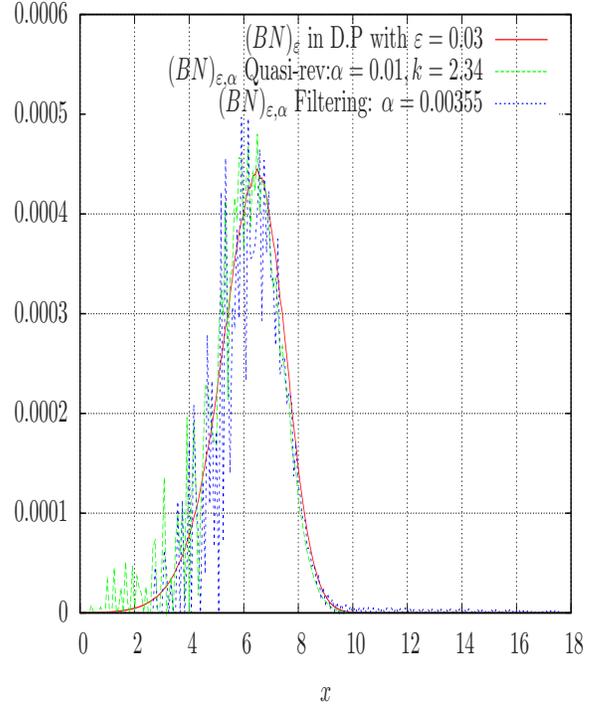}}\vfill
\resizebox{8.5cm}{9.75cm} {\input{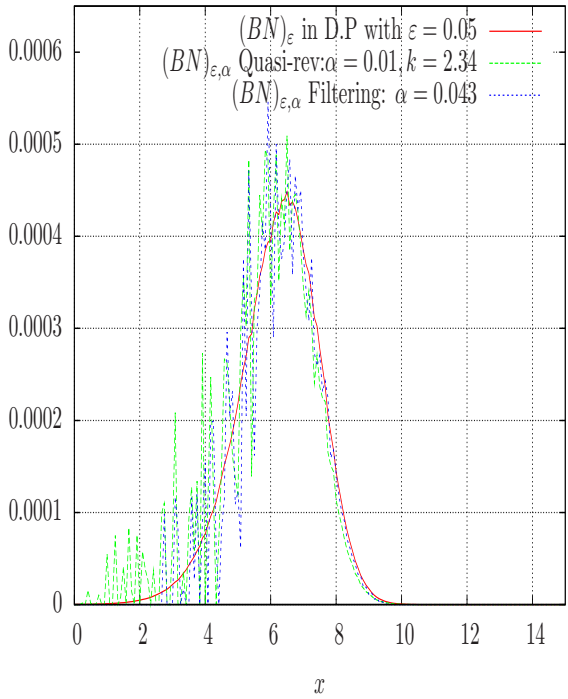}}\hfill
\resizebox{8.5cm}{9.75cm} {\input{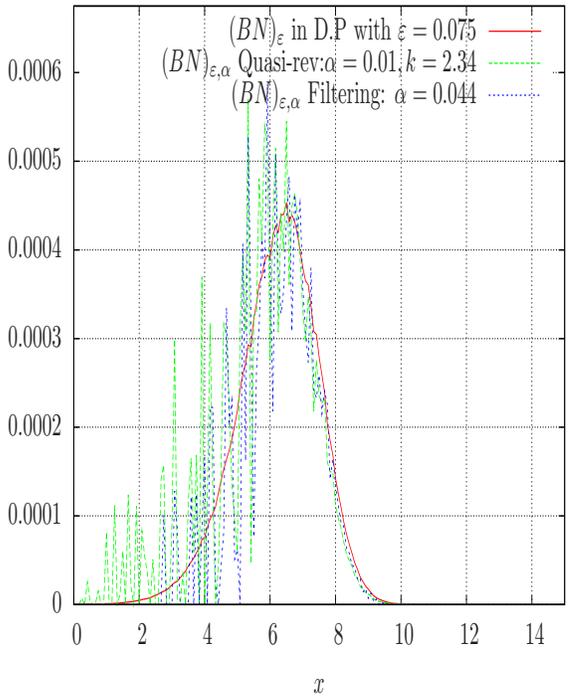}}
\caption{Numerical reconstruction of $BN$ by the measured data $N_{\varepsilon}$ for different values of $\varepsilon$ with the choice $B(x)=\exp(-0.08(x-12)^2)$, $c=0.015$ and $g(x)=x$.}
\label{fig7}

\end{figure} 
\newpage	
For various choice of the parameter $\varepsilon$ we compute the relative error thanks to the relation \eqref{eq_err} and obtain the following representations
 \begin{figure}[!ht]
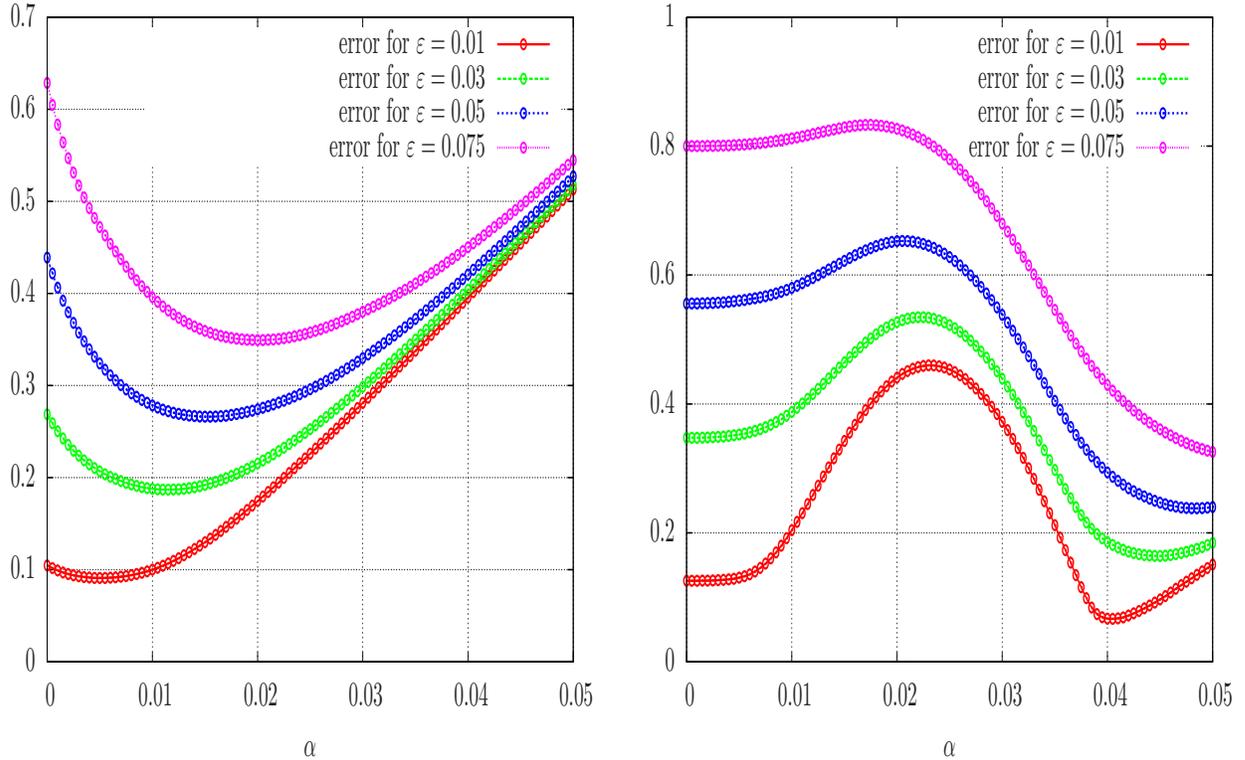

\resizebox{8.5cm}{10.5cm} {\input{leon_marie5.tex}}\hfill
\resizebox{8.5cm}{10.5cm} {\input{leon_marie6.tex}}
 \caption{Numerical errors for different values of $\varepsilon\neq0$ with $B(x)=\exp(-0.08(x-12)^2)$, $c=0.015$ and $g(x)=x$ in the direct problem.
 Left : errors by Quasi-reversibility method ($k=2.34$) . Right : errors by Filtering method.}
\label{fig8}
\end{figure} 
\\
\end{itemize}

\begin{Rmq}
Let us note that for data with high noise values $i.e.$ $\varepsilon>0.075$ the regularization by Quasi-reversibility method gives numerically better results than the Filtering one which creates big oscillations.\\ 
\end{Rmq}

\begin{itemize}
\item[\textbullet] {\bf{For the choice $\kappa(x,y)=\ds\frac{1}{y}\kappa_{0}(\frac{x}{y})$ with $\kappa_{0}\sim\mathcal{N}(\frac{1}{2},\frac{1}{4})$}} 

 \begin{figure}[!ht]
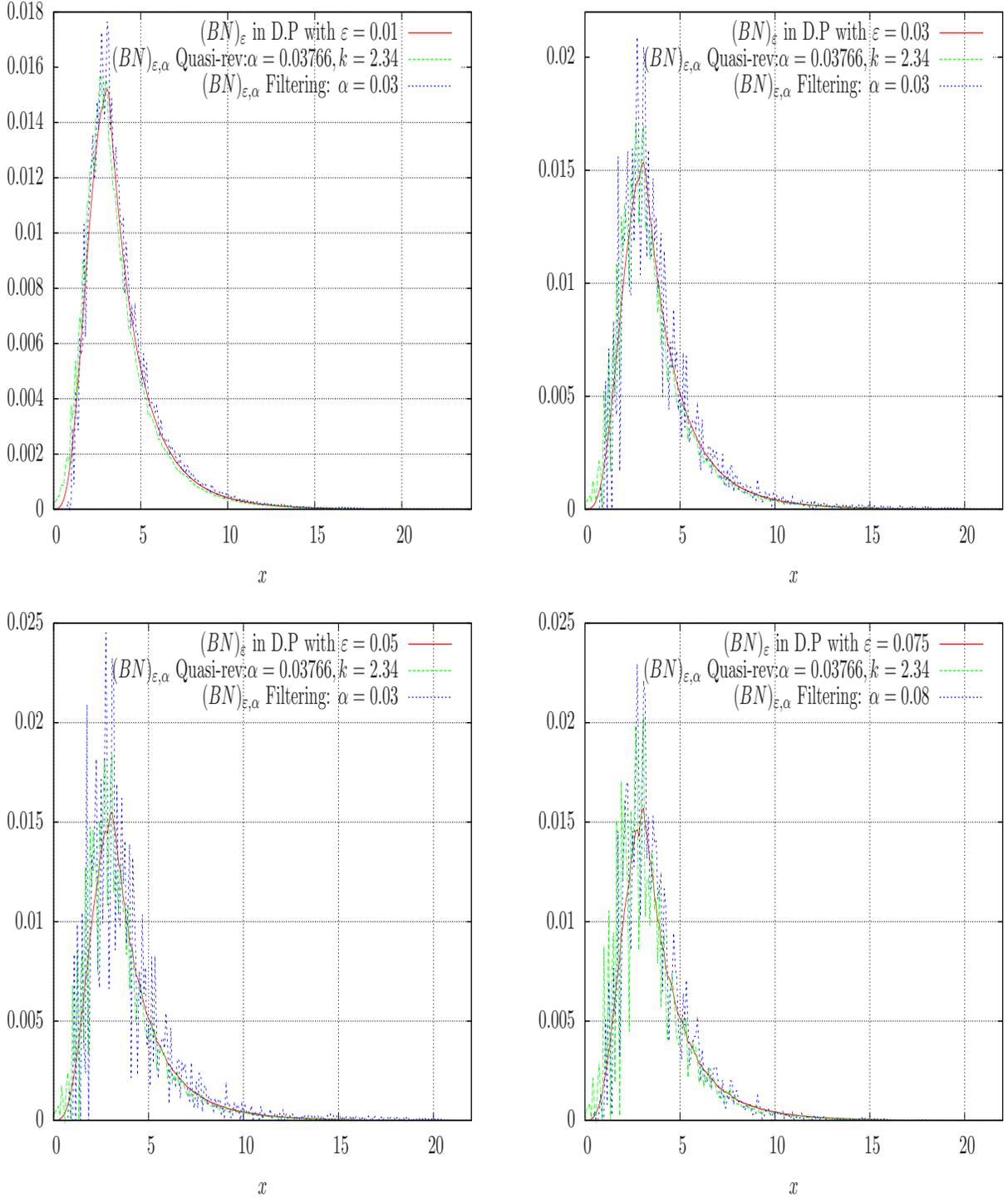

\resizebox{8.5cm}{9.75cm} {\input{comp7.tex}}\hfill
\resizebox{8.5cm}{9.75cm} {\input{comp8.tex}}\vfill
\resizebox{8.5cm}{9.75cm} {\input{comp9.tex}}\hfill
\resizebox{8.5cm}{9.75cm} {\input{comp10.tex}}
\caption{Numerical reconstruction of $BN$ by the measured data $N_{\varepsilon}$ for different values of $\varepsilon\neq0$ with the choice $B=\min(1,\frac{x^2}{10})$, $c=1$ and $g(x)=x^{1/2}$.}
\label{fig9}

\end{figure}

\begin{figure}[!ht]
\resizebox{8.5cm}{8.75cm} {\input{norme_q_epsi_second.tex}}
\resizebox{8.5cm}{8.75cm} {\input{norme_f_epsi_second.tex}}
\caption{Numerical errors for different values of $\epsilon\neq0$ with $B=\min(1,\frac{x^2}{10})$, $c=1$ and $g(x)=x^{1/2}$ in the direct problem.
 Left : errors by Quasi-reversibility method ($k=2.34$) . Right : errors by Filtering method.} 
\label{fig10}
\end{figure} 
\end{itemize}

\subsection*{Discussion}
As shown by the numerical illustrations above, and after that we tried many different shapes of regularization (trying for instance a wide variety of $k$ and $p,$ with $\pm \alpha,$ in the quasi-reversibility method), our simulations still present some delicate points. Indeed, even if the regularization methods prove to give better result than the naive "brute force" method as shown by Figure \ref{fig8}, the gain remains relatively small, and the regularizing parameter $\alpha$ has also to remain small to avoid wrong reconstructions. Due to this small regularization, as shown by Figures \ref{fig7}, \ref{fig9}, the noise is filtered but not as much as we hoped first - especially for smaller $x,$ that are farer from the departing point of the algorithm. Finally, the parameter $\alpha$ needs to stay in a confidence interval, selected, for a given growth rate $g(x),$ from a range of simulations carried out for various plausible birth rates (see for instance Figures \ref{fig4}, \ref{fig6}, \ref{fig8}). 
     
 \section{Conclusion}
 We have addressed here the problem of recovering a birth rate $B$ of a size-structured population from measurements of the time-asymptotic profile of its density, in the general case when a given individual can give birth to two daughters of inequal sizes. Compared to the work carried out in \cite{P.Z, D.P.Z, Groh} this last assumption has raised new difficulties, the principal one being that we have no other choice than considering the equation from the "viewpoint" of the daughter cell - what implies to take into account the nonlocal integral term.
   We established theoretical estimates and built numerical methods to solve it. As shown above by our numerical illustrations however, some issues still remain to be solved, especially the behavior of the algorithm for smaller $x$ and the cancellations of oscillations (also present in \cite{D.P.Z, Groh}).

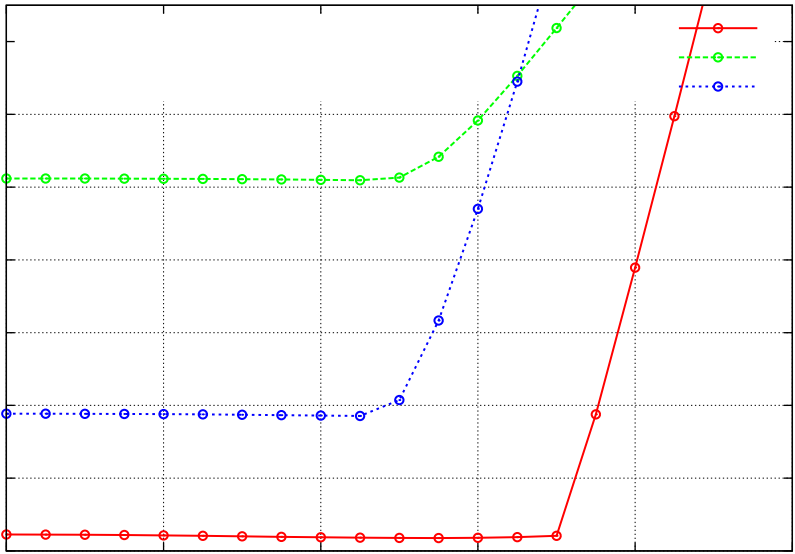
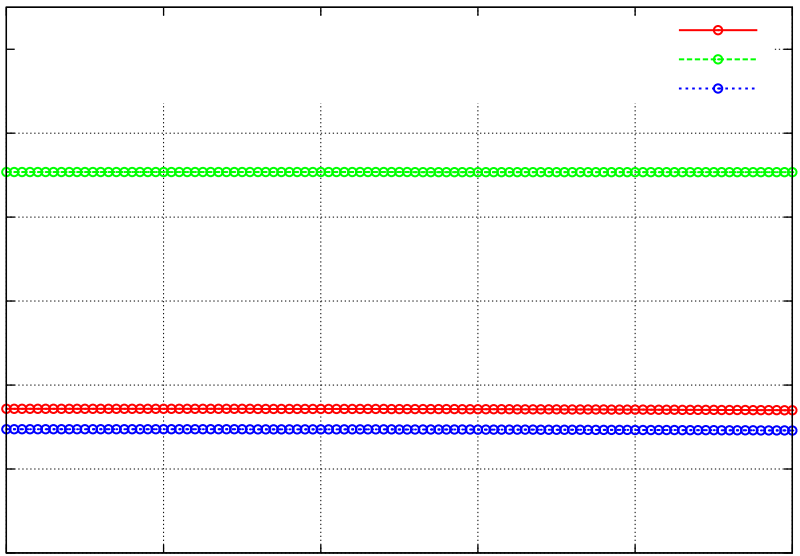
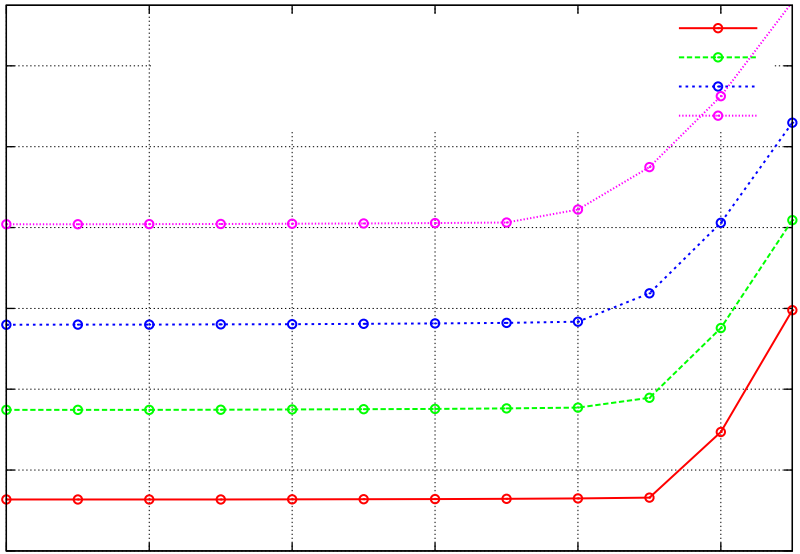
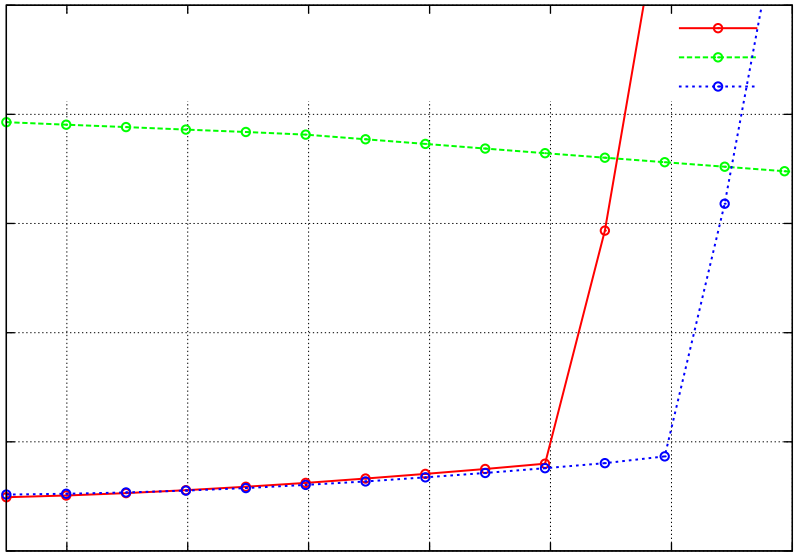
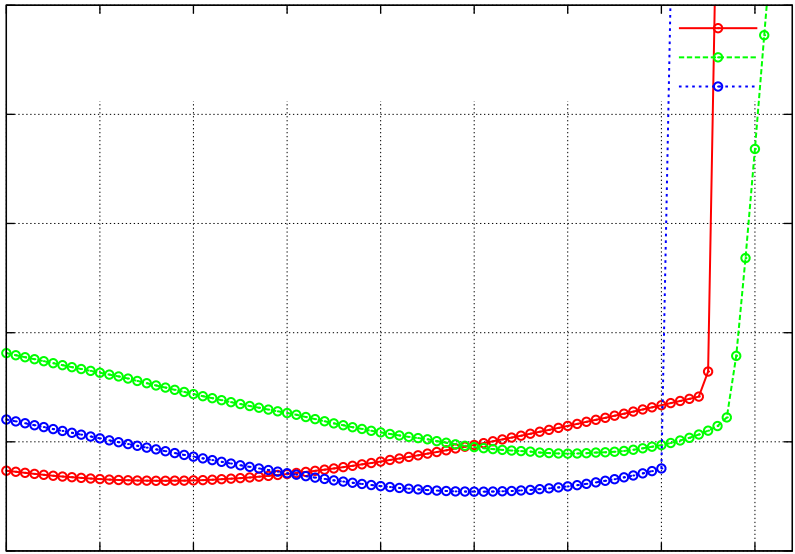
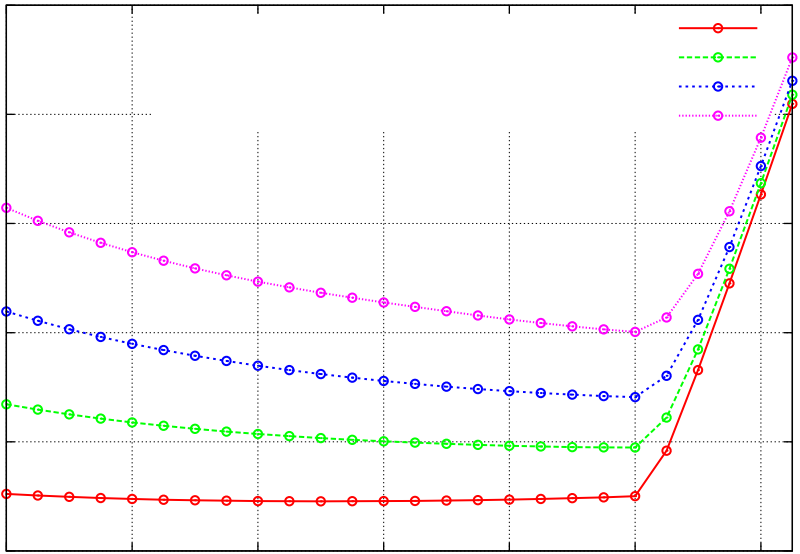
\end{document}